\documentclass[12pt,reqno]{amsart}

\usepackage{latexsym}

\newcommand{\boxb}{\raisebox{-.8pt}{$\Box$}_b}
\newcommand{\sm}{\setminus}
\newcommand{\szego}{Szeg\"o }
\newcommand{\nhat}{\raisebox{2pt}{$\wh{\ }$}}
\newcommand{\Si}{\Sigma}
\newcommand{\zetaone}{\zeta^{(1)}}
\newcommand{\zetatwo}{\zeta^{(2)}}
\newcommand{\inv}{^{-1}}
\newcommand{\kahler}{K\"ahler }
\newcommand{\sqrtn}{\sqrt{N}}
\newcommand{\wt}{\widetilde}
\newcommand{\wh}{\widehat}
\newcommand{\PP}{{\mathbb P}}
\newcommand{\R}{{\mathbb R}}
\newcommand{\C}{{\mathbb C}}

\newcommand{\CP}{\C\PP}
\renewcommand{\d}{\partial}
\newcommand{\dbar}{\bar\partial}
\newcommand{\ddbar}{\partial\dbar}

\newcommand{\D}{{\mathbf D}}
\renewcommand{\H}{{\mathbf H}}
\newcommand{\half}{{\frac{1}{2}}}
\newcommand{\vol}{{\operatorname{Vol}}}
\newcommand{\diag}{{\operatorname{diag}}}

\renewcommand{\phi}{\varphi}
\newcommand{\eqd}{\buildrel {\operatorname{def}}\over =}

\newcommand{\ccal}{\mathcal{C}}

\newcommand{\hcal}{\mathcal{H}}
\newcommand{\lcal}{\mathcal{L}}

\newcommand{\rcal}{\mathcal{R}}

\newcommand{\jcal}{\mathcal{J}}
\newcommand{\al}{\alpha}
\newcommand{\be}{\beta}
\newcommand{\ga}{\gamma}
\newcommand{\Ga}{\Gamma}
\newcommand{\La}{\Lambda}
\newcommand{\la}{\lambda}
\newcommand{\ep}{\varepsilon}
\newcommand{\de}{\delta}
\newcommand{\De}{\Delta}
\newcommand{\om}{\omega}
\newcommand{\Om}{\Omega}
\newcommand{\di}{\displaystyle}
\newtheorem{theo}{{\sc Theorem}}[section]
\newtheorem{cor}[theo]{{\sc Corollary}}
\newtheorem{lem}[theo]{{\sc Lemma}}
\newtheorem{prop}[theo]{{\sc Proposition}}

\newenvironment{rem}{\medskip\noindent{\it Remark:\/} }{\medskip}
\newenvironment{defin}{\medskip\noindent{\it Definition:\/} }{\medskip}

\title[Random almost holomorphic sections on symplectic manifolds]
{Random almost holomorphic sections of ample line bundles on symplectic
manifolds}

\author{Bernard Shiffman}
\author{Steve Zelditch}
\address{Department of Mathematics, Johns Hopkins University, Baltimore,
MD
21218, USA}
\email{shiffman@math.jhu.edu, zelditch@math.jhu.edu}

\thanks{Research partially supported by NSF grants \#DMS-9800479 (first
author) and \#DMS-9703775 (second author).}

\date{\today}

\begin{document}

\begin{abstract}  The spaces $H^0(M, L^N)$ of holomorphic sections of the 
powers of an ample line bundle $L$ over a compact K\"ahler manifold
$(M,\omega)$ have been generalized by Boutet de Monvel and Guillemin to spaces
$H^0_J(M, L^N)$ of `almost holomorphic sections' of ample line bundles over an
almost complex symplectic manifold $(M, J, \omega)$.  We consider the unit
spheres $SH^0_J(M, L^N)$ in the spaces $H^0_J(M, L^N)$, which we equip with
natural inner products.  Our purpose is to show that, in a probabilistic
sense, almost holomorphic sections behave like holomorphic sections as $N \to
\infty$.  Our first main result is that almost all sequences of sections $s_N
\in SH^0_J(M, L^N)$ are `asymptotically holomorphic' in the
Donaldson-Auroux sense that $||s_N||_{\infty}/||s_N||_{2} = O(\sqrt{\log N})$,
$||\bar{\partial} s_N||_{\infty}/||s_N||_{2} = O(\sqrt{\log N})$ and
$||\partial s_N||_{\infty}/||s_N||_{2} = O(\sqrt{N \log N})$.  Our second main
result concerns the joint probability distribution of the random variables
$s_N(z^p),\ \nabla s_N(z^p)$, $1\le p\le n$, for $n$ distinct points
$z^1,\dots, z^n$ in a neighborhood of a point $P_0\in M$.  We show that this
joint distribution has a universal scaling limit about $P_0$ as $N \to
\infty$. In particular, the limit is precisely the same as in the complex
holomorphic case.  Our methods involve near-diagonal scaling asymptotics of
the \szego projector $\Pi_N$ onto $H^0_J(M, L^N)$, which also yields proofs of
symplectic analogues of the Kodaira embedding theorem and Tian asymptotic
isometry theorem.
\end{abstract}

\maketitle

\section*{Introduction}

This paper is concerned with asymptotically
holomorphic sections of ample line bundles
over almost-complex symplectic  manifolds $(M, J, \omega).$  Such line
bundles and sections are symplectic
analogues of the usual objects in complex algebraic geometry.  Interest in
their properties has grown in
recent years because of their use by Donaldson \cite{DON.1, DON.2}, Auroux
\cite{A,A.2} and others \cite{A.3, Sik}
in proving symplectic analogues of standard results in complex geometry.
These results involve properties
of asymptotically holomorphic sections of high powers of the bundle,
particularly those involving their zero sets and the maps they define to
projective space.

We take up the study of asymptotically holomorphic sections from the
viewpoint of the microlocal analysis of
the $\bar{\partial}$ operator on a symplectic almost-complex manifold, and
define the class of `almost holomorphic sections' by a method due to  Boutet
de Monvel and  Guillemin \cite{Bou,BG}. Sections of powers of a complex
line bundle
$L^N \to M$ over $M$ are identified with equivariant functions $s$ on
the associated $S^1$-bundle $X$, and
the $\bar{\partial}$ operator is identified with the $\bar{\partial}_b$
operator on $X$.  In the non-integrable almost-complex symplectic case there
are in general no solutions of $\bar{\partial}_b {s} = 0$.  To define an
`almost
holomorphic section' $s$, Boutet de
Monvel and Guillemin
define  a
certain (pseudodifferential) $\bar{D}_j$-complex over $X$  \cite{BG}
\cite{BS}. The space $H^0_J(M, L^N)$ of almost holomorphic sections is then
defined as the space of sections corresponding to   solutions of $\bar{D}_0
{s} = 0$.  The operator   $\bar{D}_0$ is not uniquely or
even canonically defined, and it is difficult to explicitly write down
these almost holomorphic sections.  The importance, and we hope
usefulness, of these sections lies in the fact that they
typically have the properties of  asymptotically holomorphic
sections as defined by Donaldson and
Auroux, as we describe below.   We use the term `almost holomorphic' to
emphasize that a priori, they are distinct from `asymptotically
holomorphic' sections.

Our main results involve the `typical' behavior of almost holomorphic
sections in a probabilistic sense, as in our work with Bleher \cite{BSZ1,
BSZ2} and our prior work \cite{SZ} on holomorphic sections.  A wide variety of
measures could be envisioned here, and much of what we do is independent of
the precise choice of measure. However, the simplest measures are the Haar
measures on the unit spheres in the spaces $H^0_J(M, L^N)$.  To be precise, we
use a hermitian metric $h$ on $L$ and the volume form $dV =
\frac{\omega^m}{m!}$ on $M$ to endow $H^0_J(M, L^N)$ with an $\lcal^2$ inner
product. We denote by $SH^0_J(M, L^N)$ the elements of unit norm in $H_J^0(M,
L^N)$ and by $\nu_N$ the Haar probability measure on the sphere $SH^0_J(M,
L^N)$.  We also consider the essentially equivalent Gaussian measures on
$H^0_J(M, L^N)$. The theme of our work is to obtain results about almost
holomorphic sections by calculating asymptotically (as $N \to \infty$) the
probabilities that sections $s_N \in H^0_J(M, L^N)$ do various things.  This
theme has a variety of potential applications in geometry, which we hope to
pursue in the future.

In this article we focus on two applications. Our first main result
gives estimates on various norms of a typical sequence of almost holomorphic
sections of growing degree.
Let us recall that a sequence of sections
$s_N$ is called asymptotically holomorphic by Donaldson and Auroux
\cite{DON.1,A} if
$$ \|s_N\|_\infty + \|\bar{\partial} s_N\|_\infty = O(1),\  \|\nabla
s_N\|_\infty + \|\nabla \bar{\partial}
s_N|\|_\infty =
O(\sqrtn),\
\|\nabla \nabla s_N\|_\infty = O(N).
$$
We will prove that almost every sequence $\{s_N\}$ of ${\lcal^2}$-normalized
($\|s_N\|_2 = 1$) almost holomorphic sections in the Boutet de Monvel
-Guillemin sense is close to being asymptotically holomorphic in the Donaldson
sense.  We also let $\nabla:\ccal^\infty(M,L^N\otimes
(T^*M)^{\otimes k})\to \ccal^\infty(M,L^N\otimes (T^*M)^{\otimes (k+1)})$
denote the connection, and we write
$\nabla^k=\nabla\circ\cdots\circ\nabla:\ccal^\infty(M,L^N)\to
\ccal^\infty(M,L^N\otimes (T^*M)^{\otimes k} )$. We also have the
decomposition $\nabla=\d +\dbar$; note that here $\dbar$ depends on the choice
of connection.

\begin{theo}\label{aah}Endow the infinite product $\Pi_{N=1}^{\infty}
SH^0_J(M,
L^N)$ with the product spherical measure $\nu_{\infty}: = \Pi_{N=
1}^{\infty}
\nu_N.$ Then $\nu_\infty$-almost every sequence $\{s_N\}$ of sections
satisfies the following estimates:
$$\begin{array}{ll} \|s_N\|_\infty = O(\sqrt{\log N}), &\quad \|\nabla^k
s_N\|_\infty = O(N^{k/2}
\sqrt{\log N}),\\ \\  \|\bar{\partial} s_N\|_\infty = O(\sqrt{\log N}),
&\quad \|\nabla^k
\bar{\partial} s_N\|_\infty = O(N^{k/2}\sqrt{\log N}),\quad (k\ge 1)\,.
\end{array}
$$
\end{theo}

\medskip
Our second
main result concerns  the joint probability  distribution
$$\D^N_{(z^1,\dots,z^n)}=D^N(x,\xi;z^1,\dots,z^n)dxd\xi$$
 of the random variables
\begin{equation}\label{randomvar}x^p=s_N(z^p),\
\xi^p=N^{-\frac{1}{2}}\nabla s_N(z^p)\qquad (1\le p\le n)\end{equation}
on $SH^0(M,L^N)$, for
$n$ distinct points $z^1,\dots, z^n\in M$. 
We prove that upon rescaling,  this joint 
probability
distribution has
a universal limit which agrees with that of the holomorphic case determined
in \cite{BSZ2}.
\begin{theo}\label{usljpd-sphere} Let $L$ be a pre-quantum line bundle over a
$2m$-dimensional compact integral symplectic manifold $(M,\om)$. Let $P_0\in
M$ and
choose complex local coordinates $\{z_j\}$ centered at
$P_0$ so that $\om|_{P_0}$ and
$g|_{P_0}$ are the usual Euclidean \kahler form and metric respectively.
Then
$$ \D^N_{(z^1/\sqrtn,\dots, z^n/\sqrtn)} \longrightarrow \D^\infty
_{(z^1,\dots,n^n)}$$   where $\D^\infty_{(z^1,\dots,z^n)}$ is a universal
Gaussian measure supported on the holomorphic 1-jets.
\end{theo}

A technically interesting novelty in the proof is
the role of the $\bar{\partial}$ operator. In the holomorphic case,
$\D^N_{(z^1,\dots,z^n)}$ is supported
on the subspace of sections satisfying $\bar{\partial} s = 0$.  In the
almost
complex case,  sections do not satisfy this equation, so
$\D^N_{(z^1,\dots,z^n)}$ is a measure on a  higher-dimensional space of
jets. However, Theorem \ref{usljpd-sphere} says that the mass in the
`$\bar{\partial}$-directions' shrinks to zero as $N \to \infty$.

An alternate statement of Theorem \ref{usljpd-sphere} involves equipping
$H^0(M, L^N)$ with a Gaussian measure, and letting $\wt\D^N_{(z^1,\dots,z^n)}$
be the corresponding joint probability distribution on $H^0_J(M, L^N)$,
which is a Gaussian measure on the complex vector space of 1-jets of
sections.  We show (Theorem \ref{usljpd}) that these Gaussian measures
$\wt\D^N$ also have the same scaling limit $\D^\infty$, so that
asymptotically the probabilities are the same as in the holomorphic case, as
established in
\cite{BSZ2}.  To be more precise, recall that a Gaussian measure on $\R^n$
is a measure of the form
$$\ga_\De = \frac{e^{-\half\langle\De\inv
x,x\rangle}}{(2\pi)^{n/2}\sqrt{\det\De}} dx\,,$$ where $\De$ is a positive
definite symmetric $n\times n$ matrix.  It is then easy to see that
$\wt\D^N_{(z^1,\dots,z^n)}=\ga_{\De^N}\,$ where $\De^N$ is the covariance
matrix of the random variables in (\ref{randomvar}).  To deal with singular
measures, we introduce in \S \ref{s-universality} generalized Gaussians whose
covariance matrices are only semi-positive definite.  A generalized Gaussian
is simply a Gaussian supported on the subspace corresponding to the positive
eigenvalues of the covariance matrix.  The main step in the proof is to show
that the covariance matrices $\De^N$ underlying $\wt\D^N$ tend in the scaling
limit to a semi-positive matrix $\De^{\infty}$.  It follows that the scaled
distributions $\wt\D^N$ tend to a generalized Gaussian $\ga_{\De^{\infty}}$
`vanishing in the $\bar{\partial}$-directions.'

In joint work with Bleher \cite{BSZ3}, we use
this result to prove universality
of scaling limits of correlations of zeros in the setting of  almost
holomorphic sections over almost-complex symplectic manifolds.
The analysis underlying Theorem \ref{usljpd-sphere}  should also be useful for 
calculating many
other kinds of probabilities in
the setting of asymptotically holomorphic sections.  For instance, we 
believe it should
be useful for proving existence results for asymptotically holomorphic 
sections satisfying
transversality conditions.

These results are based on two essential analytical results which have an
independent interest and which we believe will have future applications.  The
first is the scaling asymptotics of the \szego kernels $\Pi_N(z, w)$, i.e.
the orthogonal projections onto $H^0_J(M, L^N)$.  To be more precise, we lift
the \szego kernels  to $X$ and the asymptotics are as follows:

\begin{quote}{\it Choose local
coordinates $\{z_j\}$ centered at a point $P_0\in M$ as in Theorem
\ref{usljpd-sphere} and choose a `preferred' local frame for $L$, which
together with the coordinates on $M$ give us `Heisenberg coordinates' on $X$
(see\/ {\rm \S\ref{s-heisenberg}}). We then have
\begin{equation}\begin{array}{l}
N^{-m}\Pi_N(P_0+\frac{u}{\sqrtn},\frac{\theta}{N};
P_0+\frac{v}{\sqrtn},\frac{\phi}{N})\\ \\ \quad = \frac{1}{\pi^m}
e^{i(\theta-\phi)+u \cdot\bar{v} - \half(|u|^2 + |v|^2)}\left[1+ \sum_{r =
1}^{K} N^{-\frac{r}{2}} b_{r}(P_0,u,v)
+ N^{-\frac{K +1}{2}}
R_K(P_0,u,v,N)\right],\end{array}\label{neardiagintro}\end{equation} where
$\|R_K(P_0,u,v,N)\|_{\ccal^j(\{|u|+|v|\le \rho\})}\le
C_{K,j,\rho}$ for $j=1,2,3,\dots$.}\end{quote}

A more precise statement will be given in Theorem \ref{neardiag}.
As more or less immediate corollaries of these scaling asymptotics, we prove
symplectic analogues
of the holomorphic Kodaira embedding theorem and Tian almost-isometry
theorem \cite{Ti}; these two results have
previously been proved by Borthwick-Uribe \cite{BU.1, BU.2} using a related
microlocal approach.  The  Borthwick-Uribe proof of the almost-complex Tian
theorem was in turn motivated by a
similar proof in the
holomorphic case in \cite{Ze}.

The proof of (\ref{neardiagintro}) is based on our second analytic result: the
construction of explicit parametrices for $\Pi$ and its Fourier coefficients
$\Pi_N$. These parametrices closely resemble those of Boutet de Monvel -
Sj\"ostrand \cite{BS} in the holomorphic case.  The construction is new but
closely follows the work of Menikoff and Sj\"ostrand \cite{MS,Sj} and of
Boutet de Monvel and Guillemin \cite{Bou,BG}.  For the sake of completeness,
we will give a fairly detailed exposition of the construction of the zeroth
term of the $\bar{D}_j$ complex and of the Szeg\"o kernel.

\section*{Guide for the reader}

For the readers' convenience, we provide here a brief outline of the paper.
We begin in \S \ref{s-aCR} by first describing some terminology from
symplectic geometry and then giving an outline of Boutet de Monvel and
Guillemin's construction
\cite{Bou, BG} of a complex of pseudodifferential operators, which replaces
the $\dbar_b$ complex in the symplectic setting. The zeroth term of this
complex is used to define sequences of almost holomorphic sections and
\szego projectors analogous to the integrable complex case (\S \ref
{s-dbarcomplex}). In \S \ref{s-parametrix}, we show that the \szego
projectors
$\Pi_N$ are complex Fourier integral operators of the same type as in the
holomorphic case, and we use this formulation to obtain the scaling
asymptotics of $\Pi_N(P_0+\frac{u}{\sqrtn},\frac{\theta}{N};
P_0+\frac{v}{\sqrtn},\frac{\phi}{N})$. Section \ref{s-kodaira} gives two
applications of these asymptotics: a proof of a `Kodaira
embedding theorem' (using
global almost holomorphic sections) for integral symplectic manifolds, and
a generalization of the asymptotic expansion theorem of \cite{Ze} to
symplectic manifolds. Section \ref{s-ah} uses the scaling asymptotics to
prove that sequences of almost holomorphic sections are almost surely (in
the probabilistic sense) asymptotically close to holomorphic (Theorem
\ref{aah}).  Finally, in \S
\ref{s-universality},  we determine the joint probability distributions
$\D^N,\ \wt\D^N$ and again apply the scaling asymptotics to prove Theorem
\ref{usljpd-sphere}.

The following chart shows the interdependencies of the sections:
$$\begin{array}{cccccc}
&&\S \ref{s-aCR}\\ &&\Downarrow \\ \S
\ref{s-kodaira}& \Longleftarrow &
\S \ref{s-parametrix} & \Longrightarrow &
\S \ref{s-ah}\\ && \Downarrow \\
&
& \S \ref{s-universality}\end{array}$$
We advise the reader who wishes to proceed quickly to the
applications in
\S\S
\ref{s-kodaira}--\ref{s-universality} that these sections depend only on
the scaling asymptotics of the
\szego kernel stated in Theorem \ref{neardiag} and the notation
and terminology given in \S\S \ref{s-acsymp}--\ref{s-dbarcomplex}.

\tableofcontents

\section{Circle bundles and almost CR geometry}\label{s-aCR}

We denote by $(M, \omega)$ a compact symplectic manifold such that
$[\frac{1}{\pi}\omega]$ is an integral
cohomology class. 
As is well known (cf.\ \cite[Prop.~8.3.1]{W}; see also \cite{GS}),
there exists a
hermitian line bundle $(L, h) \to M$ and a  metric connection $\nabla$
on
$L$ whose curvature $\Theta_L$ satisfies $\frac{i}{2}\Theta_L = \omega$.
We denote by $L^N$ the
$N^{\rm th}$ tensor power of $L$. The
`quantization' of $(M, \omega)$ at Planck constant $1/N$ should be a Hilbert
space of
polarized sections of $L^N$ (\cite[p. 266]{GS}). In the complex case,
polarized sections are
simply holomorphic sections.
The notion of polarized sections is problematic  in the
non-complex
symplectic setting, since the Lagrangean subbundle $T^{1,0}M$ defining the
complex
polarization is not integrable and there usually are no `holomorphic'
sections.  A subtle but compelling
replacement for the notion of polarized section has been proposed by Boutet
de Monvel and Guillemin \cite{Bou,BG}, and
it is this notion which we describe in this section.  For the asymptotic
analysis, it is best to view sections of
$L^N$ as functions on the unit circle bundle $X\subset L^*$; we shall
describe the `almost CR geometry' of  $X$ in \S \ref{s-heisenberg} below.

\subsection{Almost complex symplectic manifolds}\label{s-acsymp}

We begin by reviewing some terminology from almost complex symplectic
geometry. An almost complex
symplectic manifold is a symplectic manifold
$(M,
\omega)$ together with an almost complex structure $J$
satisfying the compatibility condition
$\om(Jv,Jw)=\om(v,w)$ and the positivity condition.
$\omega(v, Jv) > 0$. We give $M$  the Riemannian metric
$g(v,w)=\om(v,Jw)$. We   denote by
$T^{1,0}M,
$ resp.\
$T^{0, 1}M$,  the holomorphic, resp.\ anti-holomorphic, sub-bundle of the
complex tangent bundle $TM$;  i.e., $J = i$ on $T^{1,0}M$ and $J = -i$ on
$T^{0,1}M$.  We give $M$
local coordinates
$(x_1,y_1,\dots,x_m,y_m)$, and we write
$z_j=x_j+iy_j$. As in the integrable (i.e., holomorphic) case, we let
$\{\frac{\d}{\d z_j},
\frac{\d}{\d\bar z_j}\}$ denote the dual frame to $\{dz_j, d\bar z_j\}$.
Although in our case, the coordinates
$z_j$ are not holomorphic and consequently $\frac{\d}{\d z_j}$
is generally not in $T^{1,0}M$, we nonetheless have
$$\frac{\d}{\d z_j} =
\half \frac{\d}{\d x_j}-\frac{i}{2}\frac{\d}{\d y_j}\,,\quad
\frac{\d}{\d\bar
z_j} = \half \frac{\d}{\d x_j}+\frac{i}{2}\frac{\d}{\d y_j}\,.$$
At any point
$P_0\in M$, we can choose  a local frame
$\{\bar Z_1^M,\dots, \bar Z_m^M\}$ for
$T^{0,1}M$ near $P_0$ and coordinates about $P_0$ so that
\begin{equation} \bar Z_j^M= \frac{\d}{\d\bar z_j} + \sum_{k=1}^m B_{jk}(z)
\frac{\d}{\d z_k}\,,\quad B_{jk}(P_0)=0\,,\label{localtangentframe}
\end{equation}  and hence $\d/\d z_j|_{P_0}\in T^{1,0}(M)$. 
This is one of the properties of our `preferred coordinates' defined below.

\begin{defin} Let $P_0\in M$.  A coordinate system $(z_1,\dots,z_m)$ on a
neighborhood $U$ of $P_0$ is {\it preferred\/} at $P_0$ if
$$\sum_{j=1}^m d z_j\otimes d\bar z_j =(g-i\om)|_{P_0} \,.$$
\end{defin}

In fact, the coordinates $(z_1,\dots,z_m)$ are
preferred at $P_0$ if an only if any two of the following conditions (and
hence all three) are satisfied:

\begin{enumerate}
\item[i)] $\quad\d/\d z_j|_{P_0}\in T^{1,0}(M)$, for $1\le j\le m$,
\item[ii)] $\quad\om({P_0})=\om_0$,
\item[iii)] $\quad g({P_0} )= g_0$,
\end{enumerate}
where $\omega_0$ is the standard
symplectic form and $g_0$ is the Euclidean metric:
$$\om_0=\frac{i}{2}\sum_{j=1}^m dz_j\wedge d\bar z_j =\sum_{j=1}^m
(dx_j\otimes dy_j - dy_j\otimes dx_j)\,,\quad g_0=
\sum_{j=1}^m
(dx_j\otimes dx_j + dy_j\otimes dy_j)\,.$$
(To verify this statement, note that condition (i) is equivalent to
$J(dx_j)=-dy_j$ at
$P_0$, and use $g(v,w)=\om(v,Jw)$.)
Note that by the Darboux theorem, we can choose the coordinates so that
condition (ii) is satisfied on a neighborhood of $P_0$, but this is not
necessary for our scaling results.

\subsection{The circle bundle and Heisenberg coordinates}\label{s-heisenberg}

We now let $(M, \omega, J)$ be a compact, almost complex symplectic
manifold such that
$[\frac{1}{\pi}\omega]$ is an integral
cohomology class, and we choose a
hermitian line bundle $(L, h) \to M$ and a  metric connection $\nabla$
on
$L$ with  $\frac{i}{2} \Theta_L = \omega$.
In order to simultaneously analyze sections of all positive powers $L^N$
of the line bundle $L$, we work on the associated principal
$S^1$
bundle $X \to M$, which is defined as follows: let $\pi: L^* \to M$ denote
the
dual line bundle to $L$ with dual metric $h^*$, and put $X = \{v \in L^*:
\|v\|_{h^*} =1\}$.   We let
$\alpha$ be the the connection 1-form on $X$ given by $\nabla$; we then have
$d\alpha =\pi^* \omega$, and thus $\al$ is a contact form on
$X$, i.e., $\al\wedge (d\al)^m$ is a volume form on $X$.

We let $r_{\theta}x =e^{i\theta} x$ ($x\in X$) denote the
$S^1$
action on $X$ and denote its infinitesimal generator by
$\frac{\partial}{\partial \theta}$.
A section $s$ of $L$ determines an equivariant
function
$\hat{s}$ on $L^*$ by the rule $\hat{s}(\lambda) = \left(\lambda,
s(z)
\right)$ ($\lambda \in L^*_z, z \in M$). It is clear that if $\tau \in
\C$
then $\hat{s}(z, \tau \lambda) = \tau \hat{s}$. We henceforth restrict
$\hat{s}$ to $X$ and then the equivariance property takes the form
$\hat{s}(r_{\theta} x) = e^{i \theta}\hat{s}(x)$.  Similarly, a section
$s_N$
of $L^{N}$ determines an equivariant function $\hat{s}_N$ on $X$: put
\begin{equation} \label{sNhat}\hat{s}_N(\lambda) = \left( \lambda^{\otimes
N}, s_N(z)
\right)\,,\quad
\la\in X_z\,,\end{equation} where $\lambda^{\otimes N} = \lambda \otimes
\cdots\otimes
\lambda$;
then $\hat s_N(r_\theta x) = e^{iN\theta} \hat s_N(x)$.
We denote by $\lcal^2_N(X)$ the
space of
such equivariant functions transforming by the $N^{\rm th}$ character.

In the complex case, $X$ is a CR manifold.  In the general almost-complex
symplectic case
it is an almost CR manifold.  The {\it almost CR structure\/} is defined as
follows:
 The kernel of $\alpha$ defines a horizontal hyperplane bundle $H \subset
TX$. Using the projection $\pi: X \to M$, we may lift the splitting
$TM=T^{1,0}M
\oplus T^{0,1}M$ to a splitting $H=H^{1,0} \oplus H^{0,1}$. The almost CR
structure on $X$ is defined to be the splitting
   $TX = H^{1,0} \oplus H^{0,1} \oplus \C \frac{\partial}{\partial
\theta}$.   We  also consider a local orthonormal
frame $Z_1, \dots, Z_n$ of $H^{1,0}$ , resp.\ 
$\bar{Z}_1,
\dots,
\bar{Z}_m$ of $H^{0,1}$, and dual orthonormal coframes $\vartheta_1, \dots,
\vartheta_m,$ resp. $\bar{\vartheta}_1, \dots, \bar{\vartheta}_m$. On the
manifold $X$ we have
$d=
\d_b +\dbar_b +\frac{\partial}{\partial \theta}\otimes \alpha$, where
$\partial_b  =
\sum_{j = 1}^m {\vartheta}_j
\otimes{Z}_j$ and
 $\dbar_b  = \sum_{j = 1}^m \bar{\vartheta}_j \otimes \bar{Z}_j$.
We  define the almost-CR $\bar{\partial}_b$
operator by  $\bar{\partial}_b = df|_{H^{1,0}}$.  
Note that for an $\lcal^2$ section $s^N$ of $L^N$, we have
\begin{equation}\label{dhorizontal}
(\nabla_{L^N}s^N)\nhat = d^h\hat s^N\,,\end{equation} where
$d^h=\d_b+\dbar_b$ is the horizontal derivative on $X$.

Our near-diagonal asymptotics of the \szego kernel (\S \ref{s-neardiag})
are given in terms of the Heisenberg dilations, using local
`Heisenberg coordinates' at a point $x_0\in X$.  To describe these
coordinates, we first need the concept of a
`preferred frame':

\begin{defin}  A {\it preferred frame\/} for $L\to M$ at a point $P_0\in
M$ is a local frame $e_L$ in a neighborhood of $P_0$ such that 

\begin{enumerate}
\item[i)] $\quad \|e_L\|_{P_0} =1$;
\item[ii)] $\quad \nabla e_L|_{P_0} = 0$;
\item[iii)] $\quad \nabla^2 e_L|_{P_0} = -(g+i\om)\otimes
e_L|_{P_0}\in T^*_M\otimes T^*_M\otimes L$.
\end{enumerate} \end{defin}
(A preferred frame can be constructed by multiplying an arbitrary frame by
a function with specified 2-jet at $P_0$; any two such frames agree to
third order at $P_0$.)  Once we have property (ii), property (iii) is
independent of the choice of connection on
$T^*_M$ used to define $\nabla:\ccal^\infty(M,L\otimes T^*_M)\to
\ccal^\infty(M,L\otimes T^*_M \otimes T^*_M)$. In fact, property (iii) is a
necessary condition for obtaining universal scaling asymptotics, because of
the `parabolic' scaling in the Heisenberg group. 
 Note that if $e_L$ is a
preferred frame at $P_0$ and if $(z_1,\dots,z_m)$ are preferred
coordinates at $P_0$, then we compute the Hessian of $\|e_L\|$:
$$\left(\nabla^2
\|e_L\|_h\right)_{P_0} = \Re \left(\nabla^2 e_L,e_L\right)_{P_0}=
-g(P_0)\,;$$ thus if the preferred coordinates are `centered' at
$P_0$ (i.e., $P_0=0$), we have
\begin{equation} \label{asquared}\|e_L\|_h= 1 - \half |z|^2 +
O(|z|^3)\,.\end{equation}

\begin{rem}
Recall (\cite[\S 1.3.2]{BSZ2}) that the Bargmann-Fock representation of the
Heisenberg group acts on the space of holomorphic functions on
$(M,\om)=(\C^m,\om_0)$ that are square integrable with respect to the
weight
$h=e^{-|z^2|}$. We let $L=\C^m\times \C$ be the trivial bundle. Then the
trivializing section
$e_L(z):=(z,1)$ is a preferred frame at $P_0=0$ with respect to the
Hermitian connection $\nabla$ given by $$\nabla e_L=\d \log h \otimes e_L
= -\sum_{j=1}^m \bar z_j dz_j \otimes e_L\,.$$
Indeed, the above yields $\nabla^2e_L|_{0}= -\sum d\bar z_j\otimes dz_j
\otimes e_L(0)= -(g_0+i\om_0) \otimes e_L(0)$.\end{rem}

\medskip
The preferred frame and preferred coordinates together give us `Heisenberg
coordinates':

\begin{defin} A  {\it Heisenberg coordinate chart\/} at a point $x_0$ in
the principal bundle $X$ is a coordinate chart
$\rho:U\approx
V$ with $0\in U\subset \C^m\times \R$ and $\rho(0)=x_0\in V\subset X$ of the
form
\begin{equation}\rho(z_1,\dots,z_m,\theta)= e^{i\theta} a(z)^{-\half}
e^*_L(z)\,,\label{coordinates}\end{equation} where
$e_L$ is a preferred local frame for $L\to M$ at $P_0=\pi(x_0)$, and
$(z_1,\dots,z_m)$ are preferred
coordinates centered at $P_0$.
(Note that $P_0$ has coordinates $(0,\dots,0)$ and $e_L^*(P_0)=x_0$.)
\end{defin}

We now give some computations using local coordinates
$(z_1,\dots,z_m,\theta)$ of the form (\ref{coordinates}) for a local
frame $e_L$.  (For the moment, we do not assume they are Heisenberg
coordinates.) We write \begin{eqnarray*}a(z) &=& \|e^*_L(z)\|^2_{h^*}\ =\ 
\|e_L(z)\|^{-2}_h\,,\\
\alpha &=& d\theta + \beta\,,\qquad \be=\sum_{j=1}^m(A_jdz_j+\bar A_jd\bar
z_j)\,,\\
\nabla e_L&=& \phi \otimes e_L \,,\qquad \mbox{hence}\quad
\nabla e_L^{\otimes N}\ =\ N \phi \otimes e_L^{\otimes N}\,.
\end{eqnarray*} 

We let $\frac{\d^h}{\d z_j}\in H^{1,0}X$ denote the horizontal lift of
$\frac{\d}{\d z_j}$.
The condition $\left(\frac{\d^h}{\d z_j},\al\right) =0$
yields \begin{equation}\label{dhdzj} \frac{\d^h}{\d z_j} = \frac{\d}{\d
z_j} -A_j\frac{\d}{\d
\theta}\,,\quad \frac{\d^h}{\d\bar z_j} = \frac{\d}{\d\bar z_j}
-\bar A_j\frac{\d}{\d\theta}\,.\end{equation}  Suppose $s_N=fe_L^{\otimes
N}$ is a local section of $L^N$. Then by (\ref{sNhat}) and
(\ref{coordinates}),
\begin{equation}\label{sNhat*}\hat s_N(z,\theta) =
f(z)a(z)^{-\half}e^{iN\theta}\,.\end{equation}Differentiating (\ref{sNhat*})
and using (\ref{dhorizontal}), we conclude that
\begin{eqnarray}\phi &=& -\half d\log a -i
\beta\nonumber\\ &=& -\sum_{j=1}^m \left(\half\frac{\d\log a}{\d
z_j}+iA_j\right)dz_j
 -\sum_{j=1}^m \left(\half\frac{\d\log a}{\d \bar
z_j}+i\bar A_j\right)d\bar z_j\,.\label{connection}\end{eqnarray}

Now suppose that $(z_1,\dots,z_m,\theta)$ are Heisenberg coordinates at
$P_0$; i.e., $e_L$ is a preferred frame at $P_0$ and $(z_1,\dots,z_m)$
are preferred coordinates centered at  $P_0$ (with $P_0=0$). By
property (ii) of preferred frames, we have $\phi(0)=0$, and hence by
(\ref{connection})
\begin{equation}\label{da} da|_{0}=d\log a|_{0} =0,\end{equation}
\begin{equation}\label{Aj} A_j(0)=0\,,\quad (1\le j\le
m)\,.\end{equation}
By differentiating (\ref{connection}) and applying the
properties of preferred coordinates and frames, we further obtain
$$\sum_{j=1}^md\bar z_j\otimes dz_j = -\nabla \phi = 
\sum_{j=1}^m d \left(\half\frac{\d\log a}{\d
z_j}+iA_j\right)\otimes dz_j
+\sum_{j=1}^m d \left(\half\frac{\d\log a}{\d \bar
z_j}+i\bar A_j\right)\otimes d\bar z_j\ \mbox{at}\ 0.$$
Thus the following four equations are satisfied at $P_0=0$:
\begin{equation}\label{4equations}\begin{array}{rclrcl}\di
\half\frac{\d^2\log a}{\d z_j \d z_k} +i\frac{\d A_j}{\d z_k} &= & 0\,,
& \di\quad \half\frac{\d^2\log a}{\d z_j \d\bar z_k} +i\frac{\d A_j}{\d\bar
z_k} &= &\de^j_k\,,\\[12pt] \di
\half\frac{\d^2\log a}{\d\bar z_j \d z_k} +i\frac{\d\bar A_j}{\d z_k} &= &
0\,,&\di \quad \half\frac{\d^2\log a}{\d\bar z_j \d\bar z_k} +i\frac{\d\bar
A_j}{\d\bar z_k} &= &0\,,
\end{array}\end{equation} at $P_0$.  Solving (\ref{4equations}) and
recalling that $a(0)=1,\ da|_{0}=0$, we obtain
\begin{equation}\label{d2a}\frac{\d^2 a}{\d z_j \d z_k}(0)=0\,,\qquad
\frac{\d^2 a}{\d z_j \d\bar z_k}(0)=\de^j_k\,,\end{equation}
\begin{equation}\frac{\d A_j}{\d z_k}(0)=0\,,\qquad \frac{\d A_j}{\d\bar
z_k} = -\frac{i}{2}\de^j_k\,.\end{equation}
Hence $A_j=-\frac{i}{2}\bar z_j+O(|z|^2)$ and
\begin{equation}\frac{\d^h}{\d z_j} = \frac{\d}{\d z_j}
+\left[\frac{i}{2}\bar z_j +O(|z|^2)\right]
\frac{\d}{\d \theta}\,,\quad \frac{\d^h}{\d\bar z_j} = \frac{\d}{\d\bar z_j}
-\left[\frac{i}{2} z_j +O(|z|^2)\right]
\frac{\d}{\d \theta}
\,.\label{dhoriz}\end{equation}

\subsection{The $\bar{D}$ complex and Szeg\"o kernels}\label{s-dbarcomplex}

In the complex case, a holomorphic section $s$ of $L^N$ lifts to a
$\hat{s}\in \lcal^2_N(X)$ which satisfying $\dbar_b \hat{s} = 0.$  The
operator $\dbar_b$ extends to a complex satisfying $\dbar_b^2 = 0$,
which is a necessary
and sufficient condition for having a maximal family of CR holomorphic
coordinates.
In the non-integrable case  $\dbar_b^2 \not= 0$, and there may be no
solutions of $\dbar_b f = 0.$  To define polarized sections and their
equivariant lifts,
Boutet de Monvel \cite{Bou} and Boutet de Monvel - Guillemin \cite{BG}
defined a complex $\bar{D}_j$, which is a good replacement for
$\dbar_b$ in the non-integrable
case.  Their main result is:

\begin{theo}\label{COMPLEX} {\rm (see \cite{BG}, Lemma 14.11 and Theorem A
5.9)} There exists  an $S^1$-invariant  complex of first order
pseudodifferential operators $\bar{D}_j$ over $X$
$$0 \rightarrow C^{\infty}(\Lambda_b^{0,0})
\ {\buildrel \bar{D}_0 \over\to}\  C^{\infty}(\Lambda_b^{0,1})
\ {\buildrel \bar{D}_1 \over\to}\  \cdots \ {\buildrel \bar{D}_{m-1}
\over\longrightarrow}\
C^{\infty}(\Lambda_b^{0,m})\to 0\,,$$ where $\Lambda_b^{0,j}=\La^j
(H^{0,1}X)^*$, such that:

\begin{enumerate}
\item[{\rm i)}]  $\sigma(\bar{D}_j) =
\sigma(\bar{\partial}_b)$ to second order along $\Sigma:=\{(x,r\al_x):x\in X,
r>0\}\subset T^*X$;
\item[{\rm ii)}] The orthogonal projector  $\Pi : \lcal^2(X)
\to
\hcal^2(X)$ onto the
kernel of $\bar{D}_0$ is a complex Fourier integral operator
which is microlocally equivalent to the Cauchy-Szeg\"o projector of the
holomorphic case; 
\item[{\rm iii)}] $(\bar{D}_0,
\frac{\partial}{\partial \theta})$ is jointly elliptic. 
\end{enumerate}\end{theo}

The results stated here use only the $\bar D_0$ term of the complex; its
kernel consists of the spaces of almost holomorphic sections of the powers
$L^N$ of the line bundle $L$, as explained below. The complex $\bar D_j$ was
used by Boutet de Monvel -Guillemin 
\cite[Lemma~14.14]{BG} to show that the dimension of $H^0_J(M, L^N)$ or $
\hcal^2_N(X)$  is given by the Riemann-Roch formula (for
$N$ sufficiently
large). For our results, we need only the leading term of Riemann-Roch,
which we obtain as a consequence of Theorem \ref{tyz}(a). 
(The reader should be warned that the symbol is
described incorrectly in Lemma 14.11 of \cite{BG}.  However, it is correctly
described
in Theorem 5.9 of the Appendix to \cite{BG} and also in \cite{GU}).

We refer to the kernel $\hcal^2(X)=\ker \bar D_0
\cap \lcal^2(X)$ as the   Hardy space
 of square-integrable `almost CR functions' on $X$. The
$\lcal^2$ norm is with respect to the inner product
\begin{equation}\label{inner} \langle  F_1, F_2\rangle
=\frac{1}{2\pi}\int_X
F_1\overline{F_2}dV_X\,,\quad F_1,F_2\in\lcal^2(X)\,,\end{equation}
where \begin{equation}\label{dvx}dV_X=\frac{1}{m!}\al\wedge
(d\al)^m=\al\wedge\pi^*dV_M\,.\end{equation}

The $S^1$ action on $X$ commutes
with $\bar{D}_0$; hence $\hcal^2(X) = \bigoplus_{N
=0}^{\infty} \hcal^2_N(X)$ where $\hcal^2_N(X) =
\{ F \in \hcal^2(X): F(r_{\theta}x)
= e^{i
N \theta} F(x) \}$. We denote by $H^0_J(M, L^{ N})$ the space of sections
which
corresponds to $\hcal^2_N(X)$  under
the  map
$s\mapsto
\hat{s}$.  Elements of $H^0_J(M, L^{ N})$ are the {\it almost holomorphic
sections\/} of $L^N$. (Note that products of almost holomorphic sections are
not necessarily almost holomorphic.) We henceforth write $\hat s =s$ and
identify $H^0_J(M, L^{ N})$ with $\hcal^2_N(X)$.  Since $(\bar{D}_0,
\frac{\partial}{\partial \theta})$
is a jointly elliptic system, elements of $H^0_J(M, L^{ N})$ and
$\hcal^2_N(X)$ are smooth.  In many other respects, $H^0_J(M, L^N)$ is
analogous to the space of
holomorphic sections in the complex case.  Subsequent results will bear
this out.

We let $\Pi_N : \lcal^2(X) \rightarrow \hcal^2_N(X)$ denote the
orthogonal
projection.  The  level $N$ Szeg\"o kernel $\Pi_N(x,y)$ is defined by
\begin{equation} \Pi_N F(x) = \int_X \Pi_N(x,y) F(y) dV_X (y)\,,
\quad F\in\lcal^2(X)\,.
\end{equation} It can be given as
\begin{equation}\label{szego}\Pi_N(x,y)=\sum_{j=1}^{d_N}
S_j^N(x)\overline{ S_j^N(y)}\,,\end{equation} where
$S_1^N,\dots,S_{d_N}^N$ form an orthonormal basis of
$\hcal^2_N(X)$.

\subsection{Construction of the  Szeg\"o kernels}

In this section, we will sketch the construction of the operator $\bar D_0$
of Theorem
\ref{COMPLEX} in the
special setting of almost complex manifolds, and in so doing we will
describe the symbol of the complex in more detail. This will require the
introduction of many
objects from symplectic geometry and from the microlocal analysis of
$\bar{\partial}_b$.
We will need this material later on in the construction of a
parametrix for the Szeg\"o kernel.

\subsubsection{The characteristic variety of $\bar{\partial}_b$ }

In general, we denote by $\sigma_A$ the principal symbol of a
pseudodifferential operator $A$. To describe the principal symbol of
$\bar{\partial}_b$, we   introduce convenient local coordinates and frames.
Recalling that $H X=H^{1,0}X\oplus H^{0,1}X$, we again consider
local orthonormal frames $Z_1, \dots, Z_n$ of $H^{1,0}X$, resp.\
$\bar{Z}_1,
\dots, \bar{Z}_m$ of $H^{0,1}X$, and dual orthonormal coframes $\vartheta_1,
\dots,
\vartheta_m,$ resp. $\bar{\vartheta}_1, \dots, \bar{\vartheta}_m.$ Then we
have $\dbar_b  = \sum_{j = 1}^m \bar{\vartheta}_j \otimes \bar{Z}_j$.  Let
us define complex-valued functions on $T^*X$ by:
$$p_j(x,  \xi) = \langle Z_j(x), \xi),\;\;\;\; \bar{p}_j(x,  \xi) = \langle
\bar{Z}_j(x), \xi \rangle.$$  Then
$$\sigma_{\bar{\partial}_b}(x, \xi) = \sum_{j = 1}^m p_j(x, \xi)
\epsilon(\bar{\vartheta}_j)$$
where $\epsilon$ denotes exterior multiplication.
 We note that $ \{
\bar{p}_j, \bar{p}_k \} = \langle [\bar{Z}_j, \bar{Z}_k],
\xi \rangle$.

To state results, it is convenient to introduce the operator $\boxb :=
\bar{\partial}_b^* \bar{\partial}_b =
 \sum_{j =1}^m \bar{Z}_j^* \bar{Z}_j$ where $\bar{Z}_j^* $ is the
adjoint of the vector field regarded as a linear differential operator.
To conform to the notation of \cite{BG} we also put  $ q=\sigma(\boxb) =
\sum_{j= 1}^m |\bar{p}_j|^2.$
  The
characteristic variety $\Sigma = \{q = 0\}$
of $\bar{\partial}_b$ is the same as that of $\boxb$, namely the vertical
sub-bundle of $T^*X \to M.$
It is the conic submanifold
of $T^*X$  parametrized by  the graph of the contact form,
$\Sigma = \{(x, r \alpha_x): r > 0\} \sim X \times \R^+$.   It
follows that
$\Sigma$ is a symplectic submanifold. It is the dual (real) line bundle
to the vertical subbundle $V \subset TX$, since $\alpha (X) = G(X,
\frac{\partial}{
\partial \theta}).$

\subsubsection{The positive Lagrangean ideal $I$}

To construct the $\bar{D}_j$-complex replacing the
$\bar{\partial}_b$-complex in the non-integrable case, and
to construct the Szeg\"o kernel, we will need to study a positive Lagrangean
ideal $I$ whose generators will
define the principal symbol of $\bar{D}_0$.
 For background on positive Lagrangean
ideals,
see \cite{H}.

\begin{prop}  There exists a unique
positive  Lagrangean ideal $I$  with respect to $\Sigma$ containing
$q$. That is, there exists a unique ideal $I \subset I_{\Sigma}$ (where
$I_{\Sigma}$ is the ideal of functions vanishing on $\Sigma$) satisfying:
\medskip

\begin{itemize}

\item $I$ is closed under Poisson bracket;

\item $\Sigma$ is the set of common zeros
of $f \in I$;

\item There exist local generators $\zeta_1, \dots, \zeta_m$ such that the
matrix $\big(\frac{1}{i} \{\zeta_j, \bar{\zeta}_k\}\big)$ is positive
definite
on $\Sigma$ and that $q = \sum_{j,k} \lambda_{j \bar{k}} \zeta_j
\bar{\zeta}_k$, where $\{ \lambda_{j \bar{k}}\}$
is a hermitian positive definite matrix of functions.
\end{itemize}
\end{prop}

\begin{proof}
In the holomorphic case, $I$ is generated by the linear
functions
$\zeta_j(x, \xi) =
\langle \xi, \bar{Z}_j \rangle$.  In the general almost complex (or rather
almost CR) setting,
these functions do not Poisson commute and have to be modified. Since the
deviation of
an almost complex structure from being integrable (i.e. a true complex
structure) is measured
by the Nijenhuis bracket, it is not surprising that the generators $\zeta_j$
can be constructed
from the linear functions $\langle \xi, \bar{Z}_j \rangle$ and from the
Nijenhuis tensor. We
now explain how to do this, basically following the method of \cite{BG}.

  As a first approximation to the $\zeta_j$ we begin with the linear
functions $\zeta_j^{(1)}=\bar p_j$ on $T^*X$.  As mentioned above, the
$\zeta_j^{(1)}$ do not
generate a Lagrangean ideal in the non-integrable almost complex case,
indeed
\begin{equation}\label{p} \{\zeta_j^{(1)} , \zeta_k^{(1)} \} =
\langle \xi, [\bar{Z}_j(x), \bar{Z}_k(x)]\rangle\,.\end{equation}
However we do have that
$$\{\zeta_j^{(1)} , \zeta_k^{(1)} \} =
\{\langle \xi, \bar{Z}_j(x)\rangle, \langle \xi, \bar{Z}_k(x)\rangle \} = 0
\;\mbox{on}\; \Sigma.$$ Indeed, for $ (x, \xi) \in \Sigma$, we have $\xi = r
\alpha_x$ for some $r > 0$ so that
\begin{equation}\label{p0}\begin{array}{l}\{\langle \xi,
\bar{Z}_j(x)\rangle, \langle \xi, \bar{Z}_k(x)\rangle\} = r \alpha_x
([\bar{Z}_j(x), \bar{Z}_k(x)]) \\[10pt]\quad= r d \alpha_x (\bar{Z}_j(x),
\bar{Z}_k(x)) = r \pi^* \omega(\bar{Z}_j(x), \bar{Z}_k(x))= 0 \end{array}
\end{equation}
since $\{\bar{Z}_j\}$ forms a Lagrangean subspace for the horizontal
symplectic form $\pi^* \omega$.  Here, $\pi: X \to M$ is the natural
projection.  Moreover if we choose the local horizontal vector fields $Z_j$
to
be orthonormal relative to $\pi^* \omega$, then we also have:
\begin{equation}\label{convex}\begin{array}{l}\{\zeta_j^{(1)} ,
\bar\zeta_k^{(1)} \}(x,\xi) =
\langle \xi, [\bar{Z}_j(x), {Z}_k(x)]\rangle = r \pi^* \omega(\bar{Z}_j(x),
Z_k(x)) \\[10pt] \quad= ir\delta_j^k=i\delta_j^k p_\theta(x,\xi)\,, \qquad
(x,\xi)\in\Sigma\,.\end{array}
\end{equation}
Here, $p_\theta(x,\xi)=\langle \xi,
\frac{\d}{\d\theta}\rangle$.

Finally, we have
$$q = \sum_{j = 1}^m |\langle \xi, {Z}_j\rangle|^2= \sum_{j=1}^m
|\zetaone_j|^2\,.$$ Hence the second and third conditions on the $\zeta_j$
are satisfied by the functions $\zeta^{(1)}_j$.
Furthermore, equation (\ref{p}) tells us that the first condition is
satisfied to zero-th order for the ideal
$I_1=\big(\zeta_1^{(1)},\dots,\zeta_m^{(1)}\big)$.  In fact, let us
precisely describe the error.  We consider the orthonormal (relative to
$\om$) vector fields $Z_j^M=\pi_*Z_j$ of type (1,0) on $M$.  Recall that
the Nijenhuis tensor is given by
$$N(V,W)=\half\big([JV,JW]-[V,W]-J[V,JW]-J[JV,W]
\big)\,.$$
Hence,
\begin{equation}\label{Nijen}
N(Z_j^M,Z_k^M)=(-1-iJ)[Z_j^M,Z_k^M]=-2[Z_j^M,Z_k^M]_{(0,1)}
\eqd\sum_{p=1}^m
N^p_{jk}\bar Z_p^M\,.\end{equation}
We note that by definition, \begin{equation}\label{Nijensym}
N^p_{jk}=N^p_{kj}\,.
\end{equation} Furthermore, by the Jacobi identity
$$\{\{\zeta_j,\zeta_k\},\zeta_p\}
+\{\{\zeta_p,\zeta_j\},\zeta_k\} +\{\{\zeta_k,\zeta_p\},\zeta_j\}=0$$
applied to $(x,\al_x)\in\Si$, we have
\begin{equation}\label{Nijensym1}
N^p_{jk}+N^k_{pj}+N^j_{kp}=0\,.
\end{equation}

By (\ref{p0}) and (\ref{Nijen}), we have
\begin{equation}\label{pi} \{\zeta_j^{(1)} , \zeta_k^{(1)} \} =
\sum_{p=1}^m f_p^1\zeta_p^{(1)} +  \sum_{p=1}^m\bar N^p_{jk}
\bar\zeta_p^{(1)}\,.\end{equation}

We now argue, following \cite{BG}, that these functions can be successively
modified to satisfy the same conditions to infinite order on $\Sigma.$
The next step is to modify the functions $\zetaone_j$ by  quadratic terms so
that they satisfy the conditions $\{\zeta_j,\zeta_k\}\in I$ to first order
and the condition $q =
\sum_{j} |\zeta_j|^2$ to order 3 on
$\Sigma$. So we try to construct
new functions
$$\zetatwo_p=\zetaone_p + R_p\,,\quad R_p=\sum_{j,k}\nu_p^{jk}\bar\zetaone_j
\bar\zetaone_k$$   so that
\begin{eqnarray}\{\zetatwo_j , \zetatwo_k \} &=& \sum_{p}
f_p^2\zetatwo_p + \sum_{\al_1,\al_2}
\mu_{jk}^{\al_1\al_2}\bar\zetatwo_{\al_1}
\bar\zetatwo_{\al_2}\,;\label {i2}\\
q &=& \sum_{p}v_p^2 \zetatwo_p +\sum_\al \phi_p^\al \bar\zetatwo_{\al_1}
\bar\zetatwo_{\al_2} \bar\zetatwo_{\al_3}\bar\zetatwo_{\al_4}\,,\quad
(\al=(\al_1,\dots,\al_4))
\,.\label{ii2}\end{eqnarray}

Let us now solve (\ref{i2})--(\ref{ii2}) for the $\nu_p^{jk}$.  First of
all,
we choose $\nu_p^{jk}=\nu_p^{kj}$.  We have
$$\{\zetatwo_j,\zetatwo_k\}=\sum_{p=1}^m f_p^1\zeta_p^{(1)} +
\sum_{p=1}^m\bar
N^p_{jk}\bar\zetaone +\{\zetaone_j,R_k\} -\{\zetaone_k,R_j\} \mod
I_\Sigma^2\,.$$
By (\ref{convex}), we have
\begin{equation}\label{convex2} \{\zetaone_j,\bar\zetaone_k\}=i\de^k_j
p_\theta
\mod I_\Si\,, \end{equation} and thus
$$\{\zetaone_j,R_k\}=\sum_{p=1}^m 2i\nu_k^{pj}p_\theta\bar\zetaone_p \mod
I_\Si^2\,.$$ Therefore,
\begin{equation}\label{p2}
\{\zetatwo_j,\zetatwo_k\}= \sum_{p=1}^m f_p^1\zetatwo_p +\sum_{p=1}^m
\left(\bar N_{jk}^p +2i(\nu^{pj}_k-\nu^{pk}_j)p_\theta\right) \bar\zetaone_p
\mod I_\Si^2\,.\end{equation}
Hence $$\bar N_{jk}^p=
2i(\nu^{pk}_j-\nu^{pj}_k)p_\theta \quad \mbox{on}\ \Si\,,$$
or equivalently,
\begin{equation}\label{c2}
\nu^{pk}_j-\nu^{pj}_k =\frac{i}{2p_\theta}\bar N_{jk}^p \mod I_\Sigma\,.
\end{equation}

On the other hand, \begin{eqnarray*}q &=& \sum_p |\zetatwo_p-R_p|^2 \ =\
\sum_pv_p^2\zetatwo_p -R_p\bar\zetatwo_p\\
&=& \sum_pv_p^2\zetatwo_p  -\sum_{j,k,p} \nu_p^{jk} \bar\zetatwo_j
\bar\zetatwo_k
\bar\zetatwo_p + \sum_\al \phi_p^\al \bar\zetatwo_{\al_1}
\bar\zetatwo_{\al_2} \bar\zetatwo_{\al_3}\bar\zetatwo_{\al_4}\,.
\end{eqnarray*}
Hence (\ref{ii2}) is equivalent to
\begin{equation}\label{nusym} \nu_p^{jk}+\nu_k^{pj}+\nu_j^{kp}=0\,.
\end{equation} Using (\ref{Nijensym})--(\ref{Nijensym1}), we can solve
the equations (\ref{convex2}) and (\ref{c2}) to obtain
\begin{equation}\label{soln2} \nu_p^{jk}=\frac{i}{6p_\theta} \left(\bar
N^k_{pj} + \bar N^j_{pk}\right)\,.\end{equation}
Indeed, the solution (\ref{soln2}) is unique (modulo $I_\Sigma$) and hence
the
$R_p$ are unique modulo $I_\Sigma^3$.  In summary,
\begin{equation} \zetatwo_p=\zetaone_p + \frac{i}{3p_\theta}\sum_{j,k}
\bar N^k_{pj}\bar\zetaone_j \bar\zetaone_k\,.\label{zetatwo}\end{equation}

The passage from the $n^{\rm th}$ to the $(n+1)^{\rm st}$ step is similar,
and we refer to \cite[pp.~147--149]{BG}. \end{proof}

\begin{rem} Define $p_{\theta}(x, \xi) = \langle \xi,
\frac{\partial}{\partial \theta}\rangle.$ Since
the joint zero set of $\{\zeta_1, \dots, \zeta_m\}$ equals $\Sigma$ and
since $p_{\theta} \not= 0$ on
$\Sigma - 0$ it follows that $\{\zeta_1, \dots, \zeta_m, p_{\theta}\}$ is an
elliptic system of symbols.
\end{rem}

\subsubsection{The complex canonical relation}

Our eventual goal is to prove that $\Pi$ is a complex Fourier integral
operator and
to construct a parametrix for it. As a preliminary step we need to construct
and
describe the complex canonical relation $C$ underlying $\Pi$.  As is typical
with complex
Fourier integral operators, $C$ does not live in $T^*
X \times T^*X$ but rather in
its almost analytic extension $T^* \tilde{X} \times T^* \tilde{X}$. Here,
$\tilde{N}$
denotes the almost analytic extension of a $C^{\infty}$ manifold $N$.
Although the language of
almost analytic extensions may seem heavy, it is very helpful if one wishes
to understand
the full (complex) geometry of $C$.  When $N$ is real
analytic,  $\tilde{N}$ is the usual complexification of $X$, i.e. a complex
manifold
 in which $N$ sits as a totally real submanifold. The reader may find it
simpler to make this extra assumption. For background on almost analytic
extensions, we
refer to \cite{MeS,MS}.

Since $\pi: X \to M$ is an $S^1$ bundle over $M$, its complexification
$\tilde{\pi}: \tilde{X} \to \tilde{M}$
defines a $\C^*$ bundle over $\tilde{M}$. The connection form $\alpha$ has
an (almost) analytic continuation to
a connection $\tilde{\alpha}$ to this bundle and we may split $T \tilde{X} =
\tilde{H} \oplus \tilde{V}$,
where $\tilde{V} \to T \tilde{M}$ is the vertical subbundle of the fibration
$\tilde{X} \to \tilde{M}$ and where
$\tilde{H} \to T \tilde{M}$ is the kernel of $\tilde{\alpha}.$

The (almost) complexification of $T^* X$ is of course $T^*(\tilde{X})$. We
denote the canonical symplectic form
on $T^*X$ by $\sigma$ and that on $T^*(\tilde{X}$ by $\tilde{\sigma}$; the
notation is consistent because it is
the complexification of $\sigma.$  The symplectic cone $\Sigma$ complexifies
to $\tilde{\Sigma}$ and it remains
symplectic with respect to $\tilde{\sigma}.$ It is given by $\{(\tilde{x},
\tilde{\lambda} \tilde{\alpha}_{\tilde{x}}):
\tilde{\lambda} \in \C^*\}.$  We have a natural identification $L^* \iff
\Sigma$ given by $r x \to (x, r \alpha_x).$
We further note that the $\C^*$ bundle $L^* \to M$ is the fiberwise
complexification of the $S^1$ bundle $X \to M$,
hence $L^* \to M$ is the restriction of $\tilde{\pi}$ to $\tilde{\pi}^{-1}
(M).$ We will therefore view $L^*$ as
a submanifold of $\tilde{X}.$

\subsubsection{Definition of $C$}

Let
$\tilde{\zeta}_j$ be the almost analytic extensions of the functions
$\zeta_j$.  Then put
\begin{equation} \mathcal{ J}_+ = \{(\tilde{x}, \tilde{\xi}) \in T^*
\tilde{X}: \tilde{\zeta}_j = 0\  \forall j\}. \end{equation}
It is an involutive manifold of $T^*(\tilde{X})$ with the properties:
\begin{equation} \begin{array}{ll} (i) & (\jcal _+)_{\R} = \Sigma \\ & \\
(ii) & q |_{\jcal _+} \sim 0 \\ &  \\(iii) &
\frac{1}{i} \sigma(u, \bar{u}) > 0, \forall u \in T(\jcal _+)^{\bot}\\ &
\\(iv) &
T_{\rho}(\jcal _+) = T_{\rho} \tilde{\Sigma} \oplus \Lambda_{\rho}^+.
\end{array} \end{equation}
Here, $\Lambda_{\rho}^{\pm}$ is the sum of the eigenspaces of $F_{\rho}$,
the normal Hessian
of $q$,  corresponding to the eigenvalues
$\{ \pm i \lambda_j\}.$ The null foliation of $\mathcal{J}_+$ is given by
the joint Hamilton flow of the $\tilde{\zeta}_j$'s.

The following proposition,  proved in \cite{MS} and in
(\cite{BG}),Appendix, Lemma 4.5) defines the complex canonical relation $C$:

\begin{prop} There exists a unique strictly positive almost analytic
canonical relation
$C$
satisfying
\begin{equation} \label{C} \diag(\Sigma) \subset C \subset \jcal _+ \times
\overline{\jcal _+}.\end{equation} \end{prop}

Indeed,
\begin{equation} C = \{(\tilde{x}, \tilde{\xi}, \tilde{y}, \tilde{\eta}) \in
\jcal _+ \times \overline{\jcal _+}: (\tilde{x}, \tilde{\xi}) \sim
(\tilde{y}, \tilde{\eta})\}, \end{equation}
where $\sim$ is the equivalence relation of `belonging to the same leaf of
the null foliation of $\jcal _+.$  Thus, $C$ is the flow-out of its real
points, $\diag(\Sigma)$, under the joint Hamilton
flow of the $\tilde{\zeta}_j$'s. It is clear from the description that $C
\circ C = C^* = C,$ i.e. that $C$
is an idempotent canonical relation. It follows that $I^*(X \times X, C)$ is
a $*$-algebra.

\subsubsection{Definition of the Szeg\"o projector}

Having constructed $C$, we define  a   Szeg\"o
projector $\Pi$ associated to $\Sigma$ and $C$ to be a self-adjoint
projection $\Pi$
in the Fourier integral operator class $I^*(X \times X, C)$ with principal
symbol
$1$ (relative to the canonical 1/2-density of $C$).

It is simple to prove the existence of such a projection (see \cite{BG},
Appendix A.4):
Since  $I^*(X \times X, C)$ is a $*$-algebra, there  exists
an element  $A \in I^*(X \times X, C)$
with $\sigma_A = 1$ on $\diag(\Sigma)$, or more precisely with $\sigma_A$
equal to a projection onto a prescribed
'vacuum state'.  The principal symbols of  $A^2 - A$ and $A - A^*$ then
vanish, so these operators are
of negative order.  It follows that    the spectrum of $A$ is
concentrated near $\{0, 1\}$. Hence there exists
 an analytic function in a neighborhood of the spectrum such that
$F(A):= \Pi$ is a true projection.  Since
$I^*(X \times X, C)$ is closed under functional calculus, this projection
lies in that algebra.

We note that $\Pi$ is far from unique; given any $\Pi$ one could set $\Pi' =
e^{i
A} \Pi e^{-i A}$ where $A$ is a pseudodifferential operator of order $-1$.
We just fix one choice in what follows.

\begin{rem}  In \cite{BG} the term Szeg\"o projector (or Toeplitz structure)
is
used for a projection operator with wave front set on $\Sigma$ which is
microlocally
equivalent to the following model case on
 $\R^{2m+ 2} \times \R^{2m}$ (\cite[Sec.~5]{Bou}, \cite{BG}).

\medskip Let us use coordinates $y \in \R^{2m+2}, t \in \R^{2m}$, let $\eta,
\tau$ be the symplectically dual coordinates and
consider the operators
$$A_j := D_{y_j}+ i y_j |D_t|,\;\;\;\;\;  \bar{A}_j = D_{y_j} - i y_j
|D_t|.$$
Here, $D_x = \frac{\partial}{i \partial x}$ and $|D_t|$ is Fourier
multiplication by $|\tau|.$ The operators
$A_j$, resp.\ $\bar{A}_j$, are what are familiarly known as  creation
operators, resp.\ annihilation, operators
in the representation theory of the Heisenberg group. The characteristic
variety of the system $\{\bar{A}_j\}$ is
given by $\Sigma^0 = \{t = \tau = 0\} \equiv \R^{2m + 2}.$ The Hardy space
is given by
$\hcal^2 = \{f: \bar{A}_j f = 0,\;\;\;\forall j\}$ and the Szeg\"o kernel is
given by the complex Fourier integral kernel
$$\Pi^0(t, y, t', y') = C_m \int_{\R^{m}} e^{i \Phi} |\tau|^{m} d
\tau,\;\;\;\; \Phi = \langle t - t', \tau\rangle
+ i |\tau| (|y|^2 + |y'|^2).$$
The positive Lagrangean ideal $I$ is generated by the symbols
$\sigma({A_j})
=\zeta_j =  \eta_j + i
|\tau| y_j$.
\end{rem}

\subsubsection{Construction of the complex }

Having defined $\Pi$, one first constructs  $\bar{D}_0$  so that $\bar{D}_0
\Pi =0$. In terms of a local frame $\bar{\vartheta}_j$ of horizontal
(0,1)-forms on $X$, we may write
\begin{equation} \bar{D}_0 f = \sum_{j = 1}^m \hat{\zeta}_j(x, D) f
\bar{\vartheta}_j. \end{equation}
The coefficient operators $\hat{\zeta}_j(x, D)$ are first order
pseudodifferential operators with
principal symbols equal to $\zeta_j$ and satisfying
$$\hat{\zeta}_j \Pi \sim 0$$
modulo smoothing operators.
That is, one  `quantizes' the $\zeta_j$'s as first order pseudodifferential
operators which annihilate $\Pi.$ Let us briefly summarize their
construction (following \cite[Appendix]{BG}).

We begin with any $S^1$-equivariant symmetric first order pseudodifferential
operator $\bar{D}_0'$ with
principal symbol equal to $\sum_{j = 1}^m \zeta_j \bar{\vartheta}_j.$  Then
$\bar{D}'_0 \Pi$ is of order
$\leq 0$ so one may  find a zeroth order
pseudodifferential
operator $Q_0$ such that $\bar{D}_0' \Pi \sim Q_0 \Pi$ (modulo smoothing
operators).   Then put:
$\bar{D}_0 =
(\bar{D}_0' - Q_0) -
(\bar{D}_0' - Q_0)\Pi$. Clearly, $\bar{D}_0 \Pi = 0$ and
$\sigma({\bar{D}_0})
= \sigma(\bar{D}_0') = \sum_{j = 1}^m \zeta_j \bar{\vartheta}_j.$
The characteristic variety of $\bar{D}_0$ is then equal to $\Sigma$.  Since
$p_{\theta}$ is the
symbol of $\frac{\partial}{\partial \theta}$ and since the system
$\{\sigma_{\bar{D}_0}, p_{\theta}\}$
has no zeros in $T^*X - 0$ it follows that $\{\bar{D}_0,
\frac{\partial}{\partial \theta}\}$ is an elliptic
system.

\begin{rem}
One can then construct the higher $\bar{D}_j$ recursively so that
$\bar{D}_j
\bar{D}_{j-1} = 0$. We refer to
\cite{BG}, Appendix \S 5, for further details. \end{rem}

\section{Parametrix for the Szeg\"o projector}\label{s-parametrix}

In \cite[Theorem~3.1]{BSZ2}, we showed that for the complex case, the
scaled Szeg\"o kernel $\Pi_N$ near the diagonal is asymptotic to the Szeg\"o
kernel $\Pi^\H_1$ of level one for the reduced Heisenberg group, given by
\begin{equation}\label{Heisenberg}\Pi^\H_1(z,\theta;w,\phi) =
\frac{1}{\pi^m} e^{i(\theta-\phi)+i\Im
(z\cdot \bar w)-\half |z-w|^2}= \frac{1}{\pi^m} e^{i(\theta-\phi)+z\cdot
\bar
w-\half(|z|^2+|w|^2)}\,.\end{equation}  The method was to apply the Boutet de
Monvel-Sj\"ostrand oscillatory integral formula
\begin{equation}\label{oscint}\Pi (x,y) =
\int_0^{\infty}
e^{i t \psi(x,y)} s(x,y,t ) dt \end{equation} arising from a
parametrix construction ({\rm \cite [Th.~1.5 and \S
2.c]{BS}}). Let us recall the construction of
$\psi(x,y)$ in the integrable complex case.
Fix a local holomorphic section $e_L$ of $L$ over $U \subset M$ and define $a
\in C^{\infty}(U)$ by   $a = |e_L|^{-2}_h$.
Since $L^*|_U \approx U \times \C$ we can define local coordinates on $L^*$
by
$(z, \lambda) \approx \lambda e_L(z)$.
Then a defining function of $X \subset L^*$ is given by $\rho(z, \lambda) =
1 - |\lambda|^2 a(z)$. Define the
function $a(z,w)$ as the almost analytic extension of $a(z)$, i.e. the
solution of $\bar{\partial}_z a = 0 =
\partial_w a, a(z,z) = a(z)$ and put $\psi(x,y) = {i}(1 - \lambda
\bar{\mu} a(z,w)).$  Then $t \psi$ is a phase for
$\Pi.$

The object of this section is to show that the universal asymptotic formula
of
\cite{BSZ2} for the near-diagonal scaled Szeg\"o kernel holds for the
symplectic case (Theorem \ref{neardiag}).  To do this, we first
show that the Boutet de Monvel-Sj\"ostrand
construction can be extended to the symplectic almost-complex case.  Indeed
we
will obtain (Theorem \ref{oscintth}) an integral formula of the form
(\ref{oscint}) for the symplectic case.  In fact, our local phase function
$\psi$ will be shown to be of the form $\psi(x,y) = {i}( 1 - \lambda
\bar{\mu} a(z,w))$, where $\overline{a(w,z)} = a(z,w)$ and hence $\psi(y,x)
=
\overline{- \psi(x,y)}$.

\subsection{Oscillatory integral for $\Pi$}

In order to obtain our integral formula, we first recall
the notion of parametrizing an almost analytic Lagrangean
$\Lambda$ by a phase function. We
assume $\phi(x, \theta)$ is a regular phase function in the sense of
(\cite[Def.~3.5]{MeS}), i.e.
that it has no critical points, is homogeneous of degree one in $\theta$,
that the differentials
$d \frac{\partial \phi}{\theta_j}$ are linearly independent over $\C$ on the
set
$$C_{\phi \R} =\{(x, \theta): d_{\theta} \phi = 0\}$$
and such that $\Im \phi \geq 0 $.   We then let $\tilde{\phi}(\tilde{x},
\tilde{\theta})$ be an almost
analytic extension, put
$$C_{\tilde{\phi}} = \{(\tilde{x}, \tilde{\theta}): d_{\tilde{\theta}}
\tilde{\phi} = 0\}$$
and define the Lagrange immersion
$$\iota_{\tilde{\phi}}: (\tilde{x}, \tilde{\theta}) \in C_{\tilde{\phi}} \to
(\tilde{x},
d_{\tilde{x}} \tilde{\phi}(\tilde{x}, \tilde{\theta})).$$
The phase $\phi$ parametrizes $\Lambda$ if $\Lambda$ is the image of this
map.

The parametrix is an explicit construction of $\Pi(x, y)$ as a complex
Lagrangean kernel.
What we wish to prove now is that $C$ can
be parametrized, exactly as in the  CR case, by a phase $\lambda \psi(x,y)$
defined
on $\R^+ \times X \times X.$  This is helpful in analyzing the scaling limit
of $\Pi_N(x,y)$.
In the following we use local coordinates $(z, \lambda)$ on $L^*$ coming
from a choice of local coordinates $z$ on $M$ and a local frame $e_L(z)$ of
$L$,
and a corresponding local trivialization $(\tilde{z}, \lambda)$ of
$\tilde{X} \to \tilde{M}$. As before, we let $a=\|e_L^*\|^2$.

\begin{theo} \label{oscintth}  Let $\Pi(x,y): \lcal^2(X) \to \hcal^2(X)$ be
the Szeg\"o kernel.
  Then there exists a unique regular phase function $i t \psi(x,y)
\in C^{\infty}(\R^+ \times X \times X)$  of positive type and a
symbol $s \in
S^{m}(X \times X \times \R^+)$ of the type $$s(x,y,t) \sim
\sum_{k=0}^{\infty} t^{m-k} s_k(x,y)$$ such that
$id_x \psi |_{x = y}= -i d_y \psi
|_{x = y} = \alpha$ and
 $$\Pi (x,y) =
\int_0^{\infty}
e^{i t \psi(x,y)} s(x,y,t ) dt. $$ Furthermore, the
almost analytic extension  $\tilde{\psi} \in C^{\infty}(\tilde{X}
\times \tilde{X})$ of $\psi$ has the form
$\tilde\psi(\tilde x,\tilde y) =
i( 1 - \lambda \bar{\mu} \tilde a(\tilde z,\tilde w))$
 with $\tilde a(z,z)=a(z)$ and $\tilde a(\tilde z,\tilde w) =
\overline{\tilde a(\tilde w,\tilde z)}$.
\end{theo}

\begin{proof}

We need to  construct a function $a(z,w)$ so that $ i t \psi$ as above
parametrizes the
canonical relation $C$, i.e that $C$ is the image of the Lagrange immersion
\begin{equation}\begin{array}{l}  \iota_{\tilde{\psi}}: C_{t\tilde{\psi}}= 
\R^+
\times \{\tilde{\psi} = 0\}\to T^*(\tilde{X} \times \tilde{X})\\ \\
(t, \tilde{x}, \tilde{y}) \mapsto (\tilde{x}, t d_{\tilde{x}} \tilde{\psi};
\tilde{y}, - t
d_{\tilde{y}} \tilde{\psi}) \end{array} \end{equation}
 Since $C$ is the unique canonical relation satisfying
$ \label{C*} \diag(\Sigma) \subset C \subset \jcal _+ \times
\overline{\jcal _+}$,
 the conditions that  $\tilde{\psi}$ parametrize $C$ are the following:
\begin{enumerate}
\item  [{\rm i)}] $\{(x,y)\in X\times X:{\psi}(x,y) = 0\} = \diag(X)$;
\item  [{\rm ii)}] $d_x \psi |_{x = y}= - d_y \psi
|_{x = y} = r \alpha$ for $x, y \in X$ and for some function $r(x) > 0$;
\item   [{\rm iii)}]   $\tilde{\zeta}_j(\tilde x, d_{\tilde x} \tilde{\psi})=
0 =
\tilde{\zeta}_j(\tilde y, d_{\tilde y} \tilde{\psi})$ on $\{\tilde{\psi} =
0\}.$
\end{enumerate}
Such a $\tilde{\psi}$ is not unique, so we  require that $r\equiv 1$ in
condition (ii), i.e., 
$$d_x \psi |_{x = y}= - d_y \psi
|_{x = y} = \alpha\,.$$

Suppose we have $\tilde\psi(\tilde x,\tilde y) =
i( 1 - \lambda \bar{\mu} \tilde a(\tilde z,\tilde w))$. We observe that
$$\tilde{\psi} = 0 \iff  \tilde{a}(\tilde{z}, \tilde{w}) = (\lambda
\bar{\mu})^{-1},$$ and hence
\begin{equation} \label{first} \begin{array}{l}   id_{\tilde{x}} \tilde{\psi}
=  \bar{\mu}
\tilde{a}(\tilde{z}, \tilde{w}) {d \lambda} +
 \lambda \bar{\mu} d_{\tilde{z}} \pi^* \tilde{a}(\tilde{z}, \tilde{w}) \\ \\
\quad = {\lambda} ^{-1}{d \lambda} +
  \tilde{a}^{-1}  d_{\tilde{z}} \tilde{a} (\tilde{z}, \tilde{w})\;\;
\iff \;\; \tilde{\psi} = 0. \end{array} \end{equation}
The conditions on $a$ are therefore:
$$\left\{ \begin{array}{ll}
 a(z,w)\lambda \bar{\mu} = 1 \iff (z, \lambda ) = (w, \mu) \in X; \\ \\
(a\inv d_z  a + \la\inv d\la)|_{\diag(X)} =- (a\inv d_w  a + \la\inv
d\la)|_{\diag(X)} =  \alpha\\ \\
\tilde{\zeta}_j\left(\tilde{z}, {\lambda},\la\inv d\la +
 \tilde{a}^{-1}  d_{\tilde{z}} \tilde{a} (\tilde{z}, \tilde{w})\right) = 0 =
\tilde{\zeta}_j\left(\tilde{w}, \tilde{\mu}, \mu\inv{d \mu} +
 \tilde{a}^{-1}  d_{\tilde{w}} \tilde{a} (\tilde{z}, \tilde{w})\right),\
\forall (z,w, \lambda, \mu)
 \end{array} \right.$$

A solution $a(z,w)$ satisfying the first condition must satisfy  $a(z,z)
|\lambda|^2 = 1$ on $X$, so that $a(z,z) |\lambda|^2$ is the
local hermitian metric on $L^*$ with
unit bundle $X$, i.e. $a(z,z)=a(z)$.

We now prove that these conditions have  a unique solution near the
diagonal. We do this by
  reducing the canonical relation $C$ by
the natural  $S^1$ symmetry. The reduced relation $C_r$ has a unique
generating
function $\log a$;  the three conditions  above on $a$ will follow
automatically from this fact.

 The $S^1$ action of $X$ lifts to $T^*X$ as the Hamiltonian flow
of the function $p_{\theta}(x, \xi):= \langle \xi, \frac{\partial}{\partial
\theta} \rangle.$ The $\zeta_j$
 are invariant under this $S^1$ action, hence
\begin{equation}\label{PC} \{p_{\theta}, \zeta_j \} = 0\;\;\forall j.
\end{equation}

Now consider the level set $\{p_{\theta} = 1\} \subset T^*X$.  Dual to the
splitting $TX = H \oplus V$ we get
a splitting $T^*X = H^* \oplus V^*$, where
$$V^*(X) = \R \alpha = H^o,\;\;\;\; H^*(X) = V^o$$
where $E^o$ denotes the annihilator of a subspace $E$, i.e. the linear
functionals which vanish on $E$. Thus,
$p_{\theta} = 0$ on the horizontal space $H^*(X)$ and $p_{\theta}(\alpha) =
1$. Since
 $p_{\theta}$ is linear on the fibers of $T^*X$,
the set $\{p_{\theta} = 1\}$ has the form $\{\alpha + h: h \in H^*(X)\}$.
We also note that $p_{\theta} (d \theta) = 1$ in the  local coordinates
$(z, \theta)$ on $X$ defined by  $\lambda = e^{i \theta}$.   Hence
$\{p_{\theta} = 1\}$  may also be
identified with  $\{d \theta + h: h \in H^*(X)\}$.

Since $\{p_{\theta} = 1\}$ is a hypersurface, its null-foliation is given by
the orbits of the Hamiltonian
flow of $p_{\theta}$, i.e. by the $S^1$ action.  We use the term `reducing
by the $S^1$-action' to mean setting
$p_{\theta} = 1$ and then dividing by this action.  The reduction of $T^*X$
is thus defined by $(T^*X)_r = p_{\theta}^{-1}(1)/S^1$.  
Since $p_{\theta}^{-1}(1)$ is an affine bundle over $X$ with fiber
isomorphic to $H^*(X)\approx T^*M,$ it is clear that $(T^*X)_r\approx T^*M$ as
vector bundles over $M$.  We can obtain
a symplectic equivalence using the local coordinates $(z, \theta)$ on $X$.
Let $(p_z, p_{\theta})$ be the corresponding
symplectically dual coordinates, so that the natural symplectic form
$\sigma_{T^*X}$ on $T^*X$ is given by $\sigma_{T^*X} =
dz \wedge d p_z + d\theta \wedge d p_{\theta}$.  The notation $p_{\theta}$
is consistent with the above. Moreover,
the natural symplectic form on $T^*M$ is given locally by $\sigma_{T^*M} =
dz \wedge d p_z.$
Now define the projection
$$\chi: p_{\theta}^{-1}(1) \to T^*M, \;\;\; \chi(z, p_z, 1, p_{\theta}) =
(z, p_z). $$
This map commutes with the $S^1$ action and hence descends to the quotient
to define a local map  over $U$, still denoted $\chi$, from $(T^*X)_r \to
T^*M$. Clearly $\chi$ is symplectic.

We now reduce the canonical relation $C$. Thus we consider the $\C^* \times
\C^*$ action on $T^*\tilde{X} \times T^*\tilde{X} - 0$ generated by
$p_{\theta}(x, \xi), p_{\theta}(y, \eta).$  The reduction of $C$ is given by
$$C_r = C \cap (p_{\theta} \times p_{\theta})^{-1}(1,1) / \C^* \times
\C^*.$$  We then use $\chi \times \chi$ to  identify $C_r$ with a
(non-homogeneous) positive canonical relation in $T^*(\tilde{M} \times
\tilde{M}).$  Thus in coordinates,
\begin{equation} C_r = \{(\tilde{z}, \tilde{p_z} , \tilde{w},
\tilde{p_w}) \in T^*(\tilde{M} \times \tilde{M}): \exists \lambda, \mu,
(\tilde{z}, \lambda,  \tilde{p_z} , 1;  \tilde{w}, \mu, \tilde{p_w}, 1) \in
C\}.
 \end{equation}

Since reduction preserves real points, it is clear that
$$\begin{array}{l} (C_r)_{\R} = C_{\R} \cap (p_{\theta} \times
p_{\theta})^{-1}(1,1) / \C^* \times
\C^* \\ \\
= \{(z, p_z , z,p_z) \in \diag(T^*(M \times M)): \exists \theta \ \mbox{such
that}\
\alpha_{z, e^{i \theta}} = d\theta + p_z\}. \end{array} $$

Let us denote by $\tilde{\zeta}_{j r}$ the reductions of the functions
$\tilde{\zeta}_j$
by the $S^1$ symmetry.  Then $\tilde{\zeta}_{j r} = 0$ on either pair of
cotangent vectors
in $C_r$.  Moreover, by the uniqueness statement on $C$ it follows that
$C_r$ is the unique
canonical relation in $T^*(\tilde{M} \times \tilde{M})$ with the given set
of real points
and in the zero set of the $\tilde{\zeta}_{j r}$'s.

We now observe that $C_r$ has, at least near the diagonal, a unique global
generating
function. This holds because the natural projection
\begin{equation} C_r \subset T^*(\tilde{M} \times \tilde{M})  \to  \tilde{M}
\times \tilde{M} \end{equation}
is a local diffeomorphism near the diagonal.  Indeed, its derivative gives a
natural isomorphism
\begin{equation} T_{\rho, \rho} C_r \approx H^* \oplus H^* \approx T(\tilde{M}
\times \tilde{M})\,. \end{equation}
Therefore, there exists a  global generating function $\log \tilde{a} \in
C^{\infty}(\tilde{M} \times \tilde{M})$
i.e.
\begin{equation} C_r =\{(\tilde{z}, d_{\tilde{z}} \log \tilde{a}, \tilde{w},
d_{\tilde{w}} \log \tilde{a}), \;
\tilde{z}, \tilde{w} \in \tilde{M}\} . \end{equation}
Since $C^* = C$ it follows that $C_r^* = C^r$ and hence that $a(w,z) =
\overline{a(z,w)}.$

Working backwards, we find that
 the function $\tilde \psi(\tilde x,\tilde y) = i(1 -
\lambda \bar{\mu} \tilde a(\tilde z,\tilde w))$ satisfies the equations
$\tilde{\zeta}_j(\tilde{x}, d_{\tilde{x}} \tilde{\psi}) =
\tilde{\zeta}_j(\tilde{y}, d_{\tilde{y}} \tilde{\psi}) = 0$ on $\tilde{\psi}
= 0.$ Therefore the Lagrange immersion
\begin{equation}\begin{array}{l}  i_{\tilde{\psi}}: C_{t\tilde{\psi}}=  \R^+
\times \{\tilde{\psi} = 0\}\to T^*(\tilde{X} \times \tilde{X})\\ \\
(t, \tilde{x}, \tilde{y}) \to (\tilde{x}, t d_{\tilde{x}} \tilde{\psi};
\tilde{y}, - t
d_{\tilde{y}} \tilde{\psi}) \end{array} \end{equation}
takes its image inside $\jcal _+ \times \overline{\jcal _+}$ and reduces to
$C_r$ under
the $S^1$-symmetry.  To conclude the proof it is only necessary to show that
the real
points of the image of $i_{\tilde{\psi}}$ equal $\diag(\Sigma).$ We know
however that these
real points reduce to $(C_r)_{\R}$ and hence that $z = w$ at real points.
But we have
$$1 = \lambda \bar{\mu} a(z,w) = e^{i (\theta - \phi)}
\frac{a(z,w)}{\sqrt{a(z)} \sqrt{a(w)}},\;\;\;
\mbox{on}\;\; \{\tilde{\psi} = 0\}$$
hence when $z = w$ we have $e^{i(\theta - \phi)} = 1$ and hence $x = y$.
Since
$d_{\tilde{x}} \tilde{\psi}(x,y)|_{x = y} = \alpha_x$, it follows that the
real points
indeed equal $\diag(\Sigma).$ Therefore $t \tilde{\psi}$ parametrizes $C$
and hence
there exists a classical symbol for which $\Pi(x,y)$ has the stated
oscillatory integral
representation.

To show that the phase is of positive type, we need to describe the
asymptotics of $a(z,w)$ near the diagonal.  Note that
in the almost-complex case, we cannot describe
$a(z,w)$ as the almost analytic extension of $a(z,z)$.
(Of course, $\tilde a(\tilde z,\tilde w)$ is the almost
analytic extension of $a(z,w)$, by definition.) For our near-diagonal
asymptotics in the nonintegrable case, we instead use the following second
order expansion of
$a$ at points on the diagonal:

\begin{lem} Suppose that $(z_1,\dots,z_m)$ are preferred coordinates and
$e_L$ is a preferred frame at a point $P_0\in M$. Then the Taylor expansion of
$a(z,w)$ at $z=w=0$ is
$$a(z,w) = 1 + z\cdot \bar w + \cdots\;.$$
\label{a2}\end{lem}

\begin{proof} To begin, we recall that $a(0,0)=a(0)= \|e_L^*(P_0)\|^2=1$. To
compute the first and second order terms, we return to the equation
\begin{equation}\label{forallzw}\zeta_j\big(z,\la, \frac{d \lambda}{\lambda}
+ d_z
\log a(z,w)\big) = 0,
\;\;\;\;
\forall (z,\la;w) \in X \times M.\end{equation}
Let us write $\zeta_j= \zeta_j^{(1)}
+ R_j^{(2)},$ where  $R_j^{(2)}$ vanishes to second order on $\Sigma$ and we
recall that $\zetaone_j(\xi)=(\bar Z_j,\xi)$. Let us also Taylor
expand $\log a$:
$$\log a = L(z,w)+Q(z,w) +\cdots\,,$$ where $L$ is linear and $Q$ is
quadratic. Since $e_L$ is a preferred frame at $P_0$, it follows from
(\ref{asquared}) that
$a(z,z)=1+|z|^2 +\cdots$ and hence 
\begin{equation}\label{logadiag} L(z,z)=0\,,\qquad
Q(z,z)=|z|^2\,.\end{equation}

Since $d_z \log a|_{z = w}  +  \frac{d \lambda}{\lambda} =
\alpha\in\Sigma$, it follows from (\ref{forallzw}) that
\begin{equation}\label{Rj2}\zetaone_j\big(z,\la,\frac{d\lambda}{\lambda} +
d_z \log a\big) = -R_j^{(2)}(z, \lambda, \frac{d \lambda}{\lambda} +
d_z
\log a) = O(|z- w|^2).\end{equation} Since $a(z,w) =\overline{a(w,z)}$, we
can write
$$L(z,w)=\sum_{j=1}^m(b_jz_j+c_j\bar z_j+\bar c_jw_j+\bar b_j\bar w_j)\,.
$$ Since the $z_j$ are preferred coordinates and
$e_L$ is a preferred frame at $P_0$, we can choose the $\bar Z_j$ so that
$\bar Z_j(0)=\frac{\d}{\d\bar z_j}$ and hence by (\ref{Rj2}),
$$0=\zetaone_j\left.\big(z,\la,\frac{d\lambda}{\lambda} + d_z \log
a\big)\right|_{z=w=0,\la=1} =
\left.\left(\frac{\d}{\d
\bar z_j}, d_z \log a\right)\right|_{(0,0)}= c_j\ \forall j\,.$$  Since
$L(z,z)=0$, we have
$b_j+\bar c_j = 0$, and hence $L=0$.

To investigate the quadratic term $Q$ in (\ref{logadiag}), we write
\begin{equation}\label{findQ}(\frac{d\la}{\la}+d_z\log a)|_{(z,w)} 
=\alpha_{z} +\sum_{j=1}^m\left[z_j U'_j+\bar z_j U''_j+ w_jV'_j + \bar
w_jV''_j\right] + O(|z|^2+|w|^2)\,,
\end{equation} where  
$$\begin{array}{lcllcl}U'_j &=& \sum_{k=1}^m\left(\frac{\d^2 Q}{\d z_j
\d z_k} dz_k + \frac{\d^2 Q}{\d z_j \d\bar z_k}d\bar z_k\right)\,,\quad &
U''_j &=& \sum_{k=1}^m \left(\frac{\d^2 Q}{\d\bar z_j
\d z_k} dz_k +  \frac{\d^2 Q}{\d\bar z_j \d\bar
z_k}d\bar z_k\right)\,,\\[10pt]
V'_j &=& \sum_{k=1}^m\left(\frac{\d^2 Q}{\d w_j
\d z_k} dz_k + \frac{\d^2 Q}{\d w_j \d\bar z_k}d\bar z_k\right)\,,\quad &
V''_j &=& \sum_{k=1}^m \left(\frac{\d^2 Q}{\d\bar w_j
\d z_k} dz_k +  \frac{\d^2 Q}{\d\bar w_j \d\bar
z_k}d\bar z_k\right)\,.
\end{array}$$

Applying $\zetaone_k$ to (\ref{findQ}) and using (\ref{Rj2}) and the fact that 
$\zetaone_k(\al_z)=0$, we have 
\begin{equation}\sum_{j=1}^m \left[z_j(\bar Z_k|_z,U'_j) + \bar z_j
(\bar Z_k|_z,U''_j) + w_j(\bar Z_k|_z,V'_j) + \bar w_j
(\bar Z_k|_z,V''_j)\right]=O(|z|^2+|w|^2)\,.\label{barZk}\end{equation}
By (\ref{localtangentframe}) and (\ref{dhdzj}),
$$ \bar Z_k|_z= \frac{\d}{\d\bar z_k} + \sum_{l=1}^m B_{kl}(z)
\frac{\d}{\d z_l} +  C_k(z)\frac{\d}{\d
\theta} \,,\quad B_{kl}(0)=0\,.$$ Hence by (\ref{barZk}),
$$\frac{\d^2 Q}{\d z_j\d\bar z_k}=\left(\frac{\d}{\d\bar z_k},U'_j\right)=
(\bar Z_k|_0,U'_j)=0\,.$$  Similarly, $\frac{\d^2 Q}{\d\bar z_j\d\bar z_k}=
\frac{\d^2 Q}{\d w_j\d\bar z_k} = \frac{\d^2 Q}{\d \bar w_j\d\bar z_k} =0$.
Thus $Q(z,w)$ has no  terms containing $\bar z_k$. 
Since $Q(z,w)=\overline{Q(w,z)}$, the quadratic function
$Q$ also has no terms containing $w_k$, so we can write
$$Q(z,w)=B(z,z)+H(z,\bar w)+\overline{B(w,w)}\,,$$ where $B$, resp.\ $H$, is a
bilinear, resp.\ hermitian, form on $\C^m$. Since $Q(z,z)=|z|^2$ (recall
(\ref{logadiag})), we conclude that $B(z,z)=0$ and
hence $Q(z,w)=H(z,\bar w)=z\cdot\bar w$.
\end{proof}

To complete the proof of Theorem \ref{oscintth}, it remains to show that
the phase is of positive type; i.e., $\Im\psi\ge 0$ 
on some neighborhood of the diagonal in
$X\times X$.
Let $x\in X$ be arbitrary and choose Heisenberg
coordinates $(z,\theta)$ at $P_0=\pi(x)$ (so that $x$ has coordinates $(0,0)$).
Recalling that $\lambda = a(z)^{-\half } e^{i \theta}$ on $X$,
we have by Lemma
\ref{a2}, 
\begin{equation*} \frac{1}{i}\psi(0,0;z, \theta) = 1 -
\frac{a(0,z)}{ \sqrt{a(z)}} e^{-i\theta}
 = (1-e^{-i\theta}) +
e^{-i\theta}\left [\half |z|^2+O(|z|^3)\right ]\;. \end{equation*}
Thus, $$\Re \left[ \frac{1}{i}\psi(0,0;z, \theta)\right]\ge 
0\quad \mbox{for\ } |\theta|<\frac{\pi}{2},\
|z|<\ep\,,$$ where
$\ep$ is independent of the point $P_0\in M$.
\end{proof}

\subsection{Scaling limit of the Szeg\"o kernel}\label{s-neardiag}

The Szeg\"o kernels $\Pi_N$ are the  Fourier coefficients of $\Pi$ defined
by:
\begin{equation}\Pi_N(x,y) = \int_0^{\infty} \int_0^{2\pi}
e^{- i
N \theta}  e^{it  \psi( r_{\theta} x,y)} s(r_{\theta} x,y,t)
d\theta dt \end{equation} where $r_{\theta}$ denotes the $S^1$ action
on $X$. Changing variables $t \mapsto N t$ gives
\begin{equation} \Pi_N(x,y) = N \int_0^{\infty} \int_0^{2\pi}
e^{ i N ( -\theta + t \psi( r_{\theta} x,y))} s(r_{\theta} x,y,
Nt) d\theta dt\,.\end{equation}

We now determine the scaling limit of the Szeg\"o kernel by the argument of
\cite{BSZ2}.  For the sake of completeness, we provide the details of the
argument and add some new details on homogeneities, which are useful in
applications.  To describe the scaling limit at a point $x_0\in X$, we 
choose a Heisenberg chart
$\rho:U,0\to X,x_0$ centered at 
$P_0=\pi(x_0)\in M$.  Recall (\S \ref{s-heisenberg}) that choosing $\rho$ is
equivalent to choosing preferred coordinates centered at $P_0$ and a
preferred local frame
$e_L$ at
$P_0$.  We then write the \szego kernel $\Pi_N$ in terms of these
coordinates:
$$\Pi_N^{P_0}(u,\theta;v,\phi)=\Pi_N(\rho(u,\theta),\rho(v,\phi))\,,$$
where the superscript $P_0$ is a reminder that we are using coordinates
centered at $P_0$.  (We remark that the function $\Pi_N^{P_0}$ depends
also on the choice of preferred coordinates and preferred frame, which we
omit from the notation.)  The first term in our asymptotic formula below
says that the $N^{\rm th}$ scaled \szego kernel looks
approximately like the  Szeg\"o kernel of level one
for the
reduced  Heisenberg group (recall (\ref{Heisenberg})):
$$\Pi_N^{P_0}(\frac{u}{\sqrtn},\frac{\theta}{N};
\frac{v}{\sqrtn},\frac{\phi}{N})\approx \Pi^\H_1(u,\theta;v,\phi)=
\frac{1}{\pi^m} e^{i(\theta-\phi)+i\Im
(u\cdot \bar v)-\half |u-v|^2}\,.$$   In the following, we shall denote the 
Taylor series of a  $\ccal^\infty$ function $f$ defined in a
neighborhood of $0\in \R^K$ by
$f \sim f_0 + f_1 + f_2 + \dots$ where $f_j$ is the homogeneous polynomial
part of degree $j$.
We also denote by $R_n ^f \sim f_{n+1} + \cdots$ the remainder term
in the Taylor expansion.

The following is our main result on the scaling asymptotics of the \szego
kernels near the diagonal.  Since the result is of independent interest, we
state our asymptotic formula in a more precise form than is needed for the 
applications in this paper.

\begin{theo} \label{neardiag}
Let $P_0\in M$ and choose a Heisenberg coordinate chart about $P_0$.
Then
$$\begin{array}{l} N^{-m}\Pi_N^{P_0}(\frac{u}{\sqrtn},\frac{\theta}{N};
\frac{v}{\sqrtn},\frac{\phi}{N})\\ \\ \qquad
= \Pi^\H_1(u,\theta;v,\phi)\left[1+ \sum_{r = 1}^{K}
N^{-r/2} b_{r}(P_0,u,v)
+ N^{-(K +1)/2} R_K(P_0,u,v,N)\right]\;,\end{array}$$
where:
\begin{itemize}

\item $ b_{r} = \sum_{\alpha=0}^{2[r/2]}
\sum_{j=0}^{[3r/2]}(\psi_2)^{\alpha}Q_{r,\al,3r-2j}
\,,$ where $Q_{r,\alpha ,d}$ is homogeneous of
degree $d$ and
\begin{equation*} \psi_2(u,v) =
u \cdot\bar{v} - \half(|u|^2 + |v|^2)\,;\end{equation*}
in particular, $b_r$ has only even
homogeneity if $r$ is even, and only odd homogeneity if $r$ is odd;\\

\item $\|R_K(P_0,u,v,N)\|_{\ccal^j(\{|u|\le \rho,\ |v|\le \rho\}}\le
C_{K,j,\rho}$ for $j\ge 0,\,\rho>0$ and $C_{K,j,\rho}$ is independent of the
point
$P_0$ and choice of coordinates.
\end{itemize}
\end{theo}

\begin{proof}
We now fix $P_0$ and consider the asymptotics of
 \begin{equation}\begin{array}{l}\displaystyle \Pi_N(
\frac{u}{\sqrt{N}}, 0;  \frac{v}{\sqrt{N}}, 0)
\\[14pt]\displaystyle
\quad\quad=
 N \int_0^{\infty} \int_0^{2\pi}
 e^{ i N ( -\theta + t\psi(  \frac{u}{\sqrt{N}}, \theta; 
\frac{v}{\sqrt{N}}, 0))} s( \frac{u}{\sqrt{N}}, \theta; 
\frac{v}{\sqrt{N}}, 0), Nt) d\theta dt \,,\end{array}\end{equation}
where $\psi$ and $s$ are the phase and symbol from Theorem \ref{oscintth}
written in terms of the Heisenberg coordinates.

On $X$ we  have $\lambda = a(z)^{-\half } e^{i \phi}$.
 So for $(x,y) = (z, \phi, w, \phi') \in X \times X$, we have by Theorem
 \ref{oscintth},
\begin{equation} \psi(z, \phi, w, \phi') = i \left[1 -
\frac{a(z,w)}{ \sqrt{a(z)} \sqrt{a(w)}} e^{i
(\phi -
\phi')}\right]\;. \end{equation} It follows that
\begin{equation}\begin{array}{l}\displaystyle\psi( 
\frac{u}{\sqrt{N}},
\theta;  \frac{v}{\sqrt{N}}, 0)\\[14pt] \displaystyle \quad\quad=
i\left[1 - \frac{a( \frac{u}{\sqrt{N}},
\frac{{v}}{\sqrt{N}})}{ \sqrt{a( \frac{u}{\sqrt{N}}, 
\frac{ u}{\sqrt{N}})} \sqrt{a( \frac{v}{\sqrt{N}}, 
\frac{ v}{\sqrt{N}})}} e^{i \theta}\right] .
\end{array}\end{equation}

We observe that the asymptotic expansion of a function $f(
\frac{u}{\sqrt{N}}, 
\frac{v}{\sqrt{N}})$ in powers of $N^{-\half}$ is just the
Taylor
expansion of $f$ at  $u = v = 0$.  By Lemma \ref{a2} and the notational
convention
established above, we have
 \begin{equation}\label{taylorh2} a(
\frac{u}{\sqrt{N}}, 
\frac{v}{\sqrt{N}}) = 1 + \frac{1}{N} u\cdot \bar v +
R_3^a ( \frac{u}{\sqrt{N}},
\frac{v}{\sqrt{N}})\,,\qquad R_3^a ( \frac{u}{\sqrt{N}},
\frac{v}{\sqrt{N}})=O(N^{-3/2})\,.\end{equation}

The entire phase
\begin{equation}\label{phase} \begin{array}{l}
  t \psi( \frac{u}{\sqrt{N}},\theta;  \frac{
v}{\sqrt{N}}, 0)  -\theta  = it \left[ 1 - \frac{a(
\frac{u}{\sqrt{N}},  \frac{
v}{\sqrt{N}})}{a( \frac{u}{\sqrt{N}},  \frac{
u}{\sqrt{N}})^{\half} a( \frac{v}{\sqrt{N}},  \frac{
v}{\sqrt{N}})^{\half}} e^{i \theta}\right] -  \theta \end{array}
\end{equation}
then has the asymptotic $N$-expansion \begin{equation}
it[ 1
- e^{i \theta}] -  \theta - \frac{it}{N}\psi_2(u,v) e^{i \theta} + t
R_3^\psi(\frac{u}{\sqrtn},\frac{v}{\sqrtn}) e^{i \theta}\,. \end{equation}

As in \cite{BSZ2}, we absorb $ (i\psi_2 + NR_3^\psi)te^{i \theta}$   into
the
amplitude,
so that    $\Pi_N^{P_0}( \frac{u}{\sqrt{N}}, 0;  \frac{v}{\sqrt{N}}, 0)$
is
an oscillatory integral with phase $$\Psi(t,\theta): = it ( 1 - e^{i
\theta})
-
\theta$$ and with amplitude
$$A(t,\theta;P_0,u,v):=Ne^{ t e^{i \theta} \psi_2(u,v)  +  it e^{i \theta} N
R_3^\psi(
\frac{u}{\sqrt{N}},\frac{
v}{\sqrt{N}})} \sum_{k = 0}^{\infty} N^{m-k} t^{m-k} s_k(
\frac{u}{\sqrt{N}},  \frac{
v}{\sqrt{N}}, \theta)\,;$$
i.e., \begin{equation}\label{phase-amplitude} \Pi_N(
\frac{u}{\sqrt{N}}, 0;  \frac{v}{\sqrt{N}}, 0) =
\int_0^\infty\int_0^{2\pi} e^{iN\Psi(t,\theta)}A(t,\theta;P_0,u,v)d\theta
dt\end{equation}

Before proceeding, it is convenient to expand $\exp\left[ it e^{i \theta} N
R_3^\psi(\frac{u}{\sqrt{N}},\frac{
v}{\sqrt{N}})\right]$ in powers of $N^{-\half}$ and to keep track of the
homogeneity in $(u,v)$ of the
coefficients.  We simplify the notation by writing $g(t, \theta) := it e^{i
\theta}$.  By definition,
$$R_3^\psi( \frac{u}{\sqrt{N}}, 
\frac{v}{\sqrt{N}}) \sim N^{-3/2} \psi_3(u,v) +  N^{-2} \psi_4(u,v) + \cdots
+ N^{-d/2} \psi_d(u,v) + \cdots.$$
 We then  have
\begin{equation}\label{EXPAND} e^{  N g R_3^\psi(
\frac{u}{\sqrt{N}}, 
\frac{
v}{\sqrt{N}})}\sim \sum_{r = 0}^{\infty} N^{-r/2}
c_{r}(u,v)\,,\end{equation} where
\begin{equation}\label{EXPAND2}
c_{r} = \sum_{\la = 1 }^{r} c_{r,r+2\la}(u,v; t,\theta)\,,\ r\ge 1\,,\quad
c_0 = c_{00}=1\,,\end{equation}
with $c_{rd}$ homogeneous of degree $d$ in
$u,v$.
(The explicit formula for $c_{rd}$ is:
$$\textstyle c_{rd} = \sum \left\{ \frac{g^n}{n!}
\Pi_{j=1}^n
\psi_{a_j}(u,v):
n\ge 1,  a_j \geq 3, \sum_{j = 1}^n a_j=d,
\sum_{j=1}^n (a_j - 2)=r
\right\}\,,\ r\ge 1\,.$$  The range of $d$ is determined by the fact
that
$d =\sum_{j = 1}^n a_j = r + 2n$ with $0 \leq n \leq r$.)

We further decompose the factor $\sum_{k = 0}^{\infty} N^{m-k} t^{m-k}
s_k(
\frac{u}{\sqrt{N}},  \frac{
v}{\sqrt{N}}, \theta)$ into the homogeneous terms
$\sum_{k, \ell = 0}^{\infty} N^{m-k- \ell/2} t^{m-k} s_{k \ell}(P_0, u,v) $
where $s_{k \ell}$ is the homogeneous term of degree $\ell$ of $s_k$.
Finally, we have
\begin{eqnarray*}A &\sim & N e^{ g N R_3^\psi( \frac{u}{\sqrt{N}}, 
\frac{
v}{\sqrt{N}}, \theta)} \sum_{k = 0}^{\infty} N^{m-k} t^{m-k} s_k(
\frac{u}{\sqrt{N}},  \frac{v}{\sqrt{N}}, \theta)\\& =&
 N^{m+1} \sum_{n = 0}^{\infty} N^{-n/2} f_{n}(u, v;t, \theta,  P_0)
\;,\end{eqnarray*} $$f_{n}
= \sum_{r + \ell + 2k  = n} c_{r} s_{k \ell}
= \sum_{k=0}^{[n/2]}  t^{m - k}\left( s_{k,n-2k}
+\sum_{r=1}^{n-2k}\sum_{\la=1}^{r} c_{r,r+2\la}  s_{k,n-2k-r}\right)
=\sum_{j=0}^{[3n/2]}f_{n,3n-2j}$$
where $f_{n,d}$ is homogeneous
of degree $d$ in $(u,v)$. (The asymptotic expansion holds in the sense of
semiclassical symbols, i.e. the remainder after summing $K$ terms
is a symbol of order $m - k - (K + 1)/2$.)

 We now  evaluate the integral for $\Pi_N$ by the method of
stationary
phase as in \cite{BSZ2}.
The phase is independent of the parameters $(u,v)$ and we have
\begin{equation} \begin{array}{l} \frac{\d}{\d t} \Psi = i (1 - e^{i \theta}
)
\\[8pt]
 \frac{\d}{\d \theta}  \Psi = t e^{i \theta} - 1 \end{array} \end{equation}
so
the critical set of the phase is the point $ \{t=1, \theta = 0\}$.
The  Hessian $\Psi''$ on the critical set equals
$$ \left( \begin{array}{ll} 0 & 1  \\
1 & i  \end{array} \right)$$
so the phase is non-degenerate and the Hessian operator $L_{\Psi}$ is given
by
$$L_{\Psi}= \langle \Psi''(1, 0)^{-1} D,D\rangle =  2
\frac{\partial^2}{\partial t
\partial \theta} - i
 \frac{\partial^2}{
\partial t^2}\,.$$  We smoothly decompose the
integral into one over $|t-1|< 1$
and one over $|t-1|>\half$. Since the only critical point of the phase
occurs
at
$t = 1, \theta = 0$, the latter is rapidly decaying in $N$ and  we may
assume the integrand to be smoothly cut off
to $|t-1|< 1$.
It follows by the stationary phase method for complex
oscillatory integrals (\cite{H}, Theorem
7.7.5) that
\begin{equation}\label{MSP}  N^{-m} \Pi_N^{P_0} ( 
\frac{u}{\sqrt{N}}, 
\frac{
v}{\sqrt{N}}, \theta)  = C\sum_{j = 0}^{J}
\sum_{n = 0}^{K}  N^{-n/2 -
j}   L_j [ e^{-ig
\psi_2} f_{n}]|_{t = 1, \theta = 0}
+ \wh R_{JK}(P_0, u, v, N),\end{equation}
where $$C= N \frac{1}{\sqrt{\rm{det} (N \Psi''(1,0)/2
\pi i)}}
 =\sqrt{-2\pi i}$$ and $L_j$ is the differential operator of order $2j$ in
$(t,\theta)$ defined by
\begin{equation} L_j f(t,\theta; P_0, u, v )  =
\sum_{\nu - \mu = j}\sum_{2
\nu \geq 3 \mu}
\frac{1}{2^{\nu} i^j \mu! \nu!}
 L_{\Psi}^{\nu} [ f(t,  \theta; P_0, u, v)  (R_3^\Psi)^{\mu}(t, \theta)]
\end{equation}
with $R_3^\Psi(t,\theta)$  the third order remainder in the Taylor
expansion
of
$\Psi$ at $(t,\theta) = (1,0).$  Also, the remainder is estimated by
\begin{equation} \label{remainder} |\wh R_{JK}(P_0, u, v, N)| \leq C'
N^{- J-\frac{K+1}{2}}
\sum_{n =0}^K\sum_{|\alpha| \leq 2J+2} \sup_{t, \theta} |D^{\alpha}_{t,
\theta} e^{-ig \psi_2} f_{n} |. \end{equation}

Since $L_{\Psi}$ is a second order operator in  $(t,
\theta)$, we see that
\begin{equation} \label{SERIES}
L_j [ e^{-ig \psi_2} f_{n}]|_{t = 1, \theta = 0}= e^{\psi_2}
\sum_{\alpha \leq 2j} (\psi_2)^{\alpha} F_{n j \alpha}\,.
\end{equation}
Therefore \begin{eqnarray}\label{MSP2}N^{-m} \Pi_N^{P_0} ( 
\frac{u}{\sqrt{N}}, 
\frac{
v}{\sqrt{N}}, \theta)  &\sim
&e^{\psi_2}\sum_{n=0}^\infty \sum_{j=0}^\infty \sum_{\al=0}^{2j}(\psi_2)^\al
N^{-\frac{n}{2}-j} F_{nj\al}\nonumber\\ &\sim
&e^{\psi_2} \sum_{r=0}^\infty
\sum_{j=0}^{[r/2]}\sum_{\al=0}^{2j}(\psi_2)^\al
N^{-r/2}F_{r-2j,j,\al}\nonumber\\&\sim
&e^{\psi_2} \sum_{r=0}^\infty \sum_{\al=0}^{2[r/2]}(\psi_2)^\al
N^{-r/2}Q_{r\al}\,.\end{eqnarray}

Thus, as with $f_{n}$ we have the homogeneous expansion:
\begin{equation} \label{Q} Q_{r\al}=\sum_{j=0}^{[3r/2]}Q_{r,\al,3r-2j}\,.
\end{equation}
Here, $Q_{r,\al,d}$ is homogeneous of degree $d$ in $(u,v)$. (The term
$\psi_2$ is distinguished by being `holomorphic' in $u$ and
`anti-holomorphic' in $v$ in a sense to be elaborated below.)
Thus we have the desired Taylor series.  The estimate for the remainder
follows from (\ref{remainder}). 
\end{proof}

\section{Kodaira embedding and Tian almost isometry
theorem}\label{s-kodaira}

\begin{defin} By the Kodaira maps we mean the maps
$\Phi_N : M \to PH^0(M,L^N)'$ defined by $\Phi_N(z) = \{s^N: s^N(z) =
0\}$. Equivalently, we can choose an orthonormal basis
$S^N_1,\dots,S^N_{d_N}$
of $H^0(M,L^N)$ and write
\begin{equation}\label{Kmap} \Phi_N : M \to\CP^{d_N-1}\,,\qquad
\Phi_N(z)=\big(S^N_1(z):\dots:S^N_{d_N}(z)\big)\,.\end{equation}
We also define the lifts of the Kodaira maps:
\begin{equation}\label{lift}\wt{\Phi}_N : X \to
\C^{d_N}\,,\qquad
\wt\Phi_N(x)=(S^N_1(x),\dots,S^N_{d_N}(x))\,.\end{equation}
\end{defin}

\medskip Note that
\begin{equation}\label{PiPhi} \Pi_N(x,y)=\wt\Phi_N(x) \cdot
\overline{\wt\Phi_N(y)}\,;\end{equation} in particular,
\begin{equation}\label{PiPhi2}
\Pi_N(x,x)=\|\wt\Phi_N(x)\|^2\,.\end{equation}

We now prove the following generalization to the symplectic category of
the asymptotic expansion
theorem of \cite{Ze} (also proved independently by \cite{Cat} using the
Bergman kernel in place of the \szego kernel) and Tian's approximate
isometry theorem
\cite{Ti}:

\begin{theo}\label{tyz} Let  $L \to (M, \omega)$ be the pre-quantum line
bundle
over a $2m$-dimensional
symplectic manifold, and let $\{\Phi_N\}$ be its Kodaira maps.  Then:

\noindent{\rm (a)} There exists a complete asymptotic expansion:
$$ \Pi_N(z,0;z,0)  =  a_0 N^m +
a_1(z) N^{m-1} + a_2(z) N^{m-2} + \dots$$
for certain smooth coefficients $a_j(z)$ with $a_0 = \pi^{-m}$.
 Hence, the  maps $\Phi_N$  are well-defined for $N\gg 0$.

\noindent{\rm (b)}
 Let $\omega_{FS}$ denote the Fubini-Study form on $\CP^{d_N-1}$.
Then $$\|\frac{1}{N}  \Phi_N^*(\omega_{FS}) -
\omega\|_{\ccal^k} = O(\frac{1}{N})$$ for any $k$.  \end{theo}

\begin{proof}
(a) Using the expansion of Theorem
\ref{neardiag} with $u=v=0$ and noting that $b_r(z,0,0)=0$ for $r$ odd,
we obtain the above expansion of $ \Pi_N(z,0;z,0)$ with
$a_r(z)=b_{2r}(z,0,0)$. (The expansion also follows by precisely
the same proof as in
\cite{Ze}.)

\smallskip\noindent (b)  In the holomorphic case, (b) followed by
differentiating (a), using
that $\Phi_N^* ( {\partial} \bar\partial \log |\xi|^2) =  {\partial}
\bar\partial
\log |\Phi_N|^2$.  In the almost complex case, $\Phi_N^*$ does
not commute with
the complex derivatives, so we need to modify the proof.  To do so, we
use the following notation: the exterior derivative on a product
manifold $Y_1\times Y_2$ can be decomposed as $d=d^1+d^2$, where
$d^1$ and
$d^2$ denote exterior differentiation on the first and second factors,
respectively. (This is formally analogous to the decomposition $d=\d +\dbar$;
e.g., $d^1d^1=d^2d^2=d^1d^2+d^2d^1=0$.) 

Recall that the Fubini-Study form $\omega_{FS}$  on
$\CP^{m-1}$ is induced by the
2-form  $\wt\omega_m= \frac{i}{2} \ddbar\log |\xi|^2$ on $\C^m\sm\{0\}$.
We consider the 2-form $\Om$ on
$(\C^m\sm\{0\})\times(\C^m\sm\{0\})$ given by
$$\Om=\frac{i}{2}\ddbar \log \zeta\cdot\bar\eta=\frac{i}{2}d^1 d^2 \log
\zeta\cdot\bar\eta\,.$$  Note that $\Om$ is smooth on a neighborhood of the
diagonal $\{\zeta=\eta\}$, and
$$\Om|_{\zeta=\eta}=\wt\om_m$$ (where the restriction to $\{\zeta=\eta\}$
means the pull-back under the map $\zeta\mapsto (\zeta,\zeta)$).

It suffices to show that
$$\frac{1}{N} \tilde{\Phi}_N^* \omega_{d_N} \to \pi^* \omega, \;\;\;\; \pi:
X \to M.$$  To do this, we consider the maps
$$\Psi_N=\wt\Phi_N\times\wt\Phi_N:X\times X \to \C^{d_N}\times
\C^{d_N}\,,\quad
\Psi_N(x,y)=  (\wt\Phi_N(x),\wt\Phi_N(y))\,.$$
It is
elementary to check that
$\Psi_N^*$ commutes with $d^1$ and $d^2$.
By (\ref{PiPhi}),
we have
$$\Psi_N^*(\log \zeta\cdot\bar\eta) = (\log \zeta\cdot\bar\eta)\circ\Psi_N
=\log\Pi_N\,.$$

Therefore, \begin{equation}\label{PsiN}\frac{1}{N} \tilde{\Psi}_N^*
\Omega_{d_N}=
\frac{i}{2N}{\Psi}_N^*d^1 d^2 \log \zeta\cdot\bar\eta
=\frac{i}{2N} d^1d^2 {\Psi}_N^* \log \zeta\cdot\bar\eta
=\frac{i}{2N} d^1d^2 \log\Pi_N\,.\end{equation}

Restricting (\ref{PsiN}) to the diagonal, we then have
$$\frac{1}{N} \tilde{\Phi}_N^* \omega_{d_N} =
\frac{i}{2N} (d^1d^2 \log\Pi_N)|_{x=y}= \diag^*(d^1d^2 \log\Pi_N)\,,$$
where $\diag:X\to X\times X$ is the diagonal map $\diag(x)=(x,x)$.

Using Heisenberg coordinates as in Theorem \ref{neardiag}, we have by the
near-diagonal scaling asymptotics
\begin{eqnarray*} \left.\frac{1}{N} \tilde{\Phi}_N^*
\omega_{d_N}\right|_{P_0}&=& \left.\frac{i}{2N} \diag^*d^1d^2
\log\Pi_N^{P_0}(\frac{u}{\sqrtn},\frac{\theta}{N};
\frac{v}{\sqrtn},\frac{\phi}{N})\right|_0\\
&=&\left.\frac{i}{2N}\diag^* d^1d^2
\log\Pi_1^\H(u,\theta;v,\phi)\right|_0 +
O(N^{-\half})\,.\end{eqnarray*} Finally,
\begin{eqnarray} \left.
\frac{i}{2N}\diag^* d^1d^2
\log\Pi_1^\H(u,\theta;v,\phi)\right|_0 &=&\frac{i}{2N}
\diag^*d^1d^2\big[ i(\theta-\phi) + u\cdot\bar v -
\half(|u|^2+|v|^2)\big]\nonumber \\
&=& \frac{i}{2N} \sum_{q=1}^{m}du_q\wedge d\bar u_q \ = \
\frac{i}{2} \sum_{q=1}^{m}dz_q\wedge d\bar z_q
\ = \ \om|_{P_0}\,.\label{dxdy}\end{eqnarray}
\end{proof}

\begin{rem}  A more explicit way to show (b) is to expand the Fubini-Study
form:
$$\wt\omega_m=
\frac{i}{2}|\xi|^{-4}\left[ |\xi|^2 \sum_{j =1}^m d\xi_j \wedge
d\bar{\xi}_j -\sum_{j, k = 1}^m \bar{\xi}_j \xi_k d\xi_j \wedge
d\bar{\xi}_k \right]\,.$$ Then
$$\frac{1}{N} \tilde{\Phi}_N^* \omega_{d_N} =
\frac{i}{2}\Pi_N(x,x)^{-2}\{ (\Pi_N(x,x) d^1 d^2 \Pi_N(x,y) -
d^1 \Pi_N(x,y) \wedge d^2 \Pi_N(x,y)\}|_{x = y}\,,$$
and (b) follows from a short computation using Theorem
\ref{neardiag} as above.
\end{rem}

\medskip
It follows from Theorem \ref{tyz}(b) that  $\Phi_N$ is an immersion for
$N\gg
0$.  Using in part an idea
of Bouche \cite{Bouche}, we give a simple proof of the `Kodaira embedding
theorem'
for symplectic manifolds:

\begin{theo}\label{kodaira} For $N$ sufficiently large, $\Phi_N$ is an
embedding.
\end{theo}

\begin{proof} Let $\{P_N, Q_N\}$ be any sequence of distinct points
such that
$\Phi_N(P_N) = \Phi_N(Q_N)$.  By passing to a subsequence we may assume that
one of the following two cases holds:

\begin{enumerate}

\item [(i)]  The distance $r_N:= r(P_N, Q_N)$ between $P_N, Q_N$ satisfies
$r_N \sqrt{N} \to \infty;$

\item [(ii)]  There exists a constant $C$ independent of $N$ such that $r_N
\leq C \sqrt{N}.$
\end{enumerate}

\smallskip
To prove that case (i) cannot occur, we observe that
$$\int_{B(P_N, r_N)} |N^{-m}\Pi_N^{P_N} |^2 dV \geq 1 - o(1)$$
where $\Pi_N^{P_N}(x) = \Pi_N(\cdot, P_N)$ is
the `peak section' at $P_N$.
The same inequality holds for $Q_N$.  If $\Phi_N(P_N) = \Phi_N(Q_N)$ then
the
total $\lcal^2$-norm of $\Pi_N(x, \cdot)$
would have to be $\sim 2 N^m$, contradicting the asymptotic  $\sim N^m$
from Theorem \ref{tyz}(a).

To prove that case (ii) cannot occur, we assume on the contrary that
$\Phi_N(P_N)=
\Phi_N(Q_N)$, where $P_N=\rho_N(0)$ and $ Q_N
=\rho_N(\frac{v_N}{\sqrtn})$, $0\ne |v_N|\le C$, using a Heisenberg
coordinate chart $\rho_N$ about
$P_N$. We consider the function
\begin{equation}\label{fN}
f_N(t)=\frac{|\Pi_N^{P_N}(0,\frac{tv_N}{\sqrtn})|^2}{\Pi_N^{P_N}(0,0)
\Pi_N^{P_N}(\frac{tv_N}{\sqrtn},\frac{tv_N}{\sqrtn})}\,.\end{equation}
Recalling that
$$\Pi_N(x,y)= \tilde\Phi_N(x)\cdot\overline{\tilde\Phi_N(y)}\,,$$ we see
that
$f_N(0)=1$, which is a global and strict local maximum of $f_N$;
furthermore, since  $\Phi_N(P_N)=
\Phi_N(Q_N)$, we also have $f_N(1)=1$. Thus for some value of $t_N$ in the
open interval $(0,1)$, we have $f_N''(t_N)=0$.  By Theorem
\ref{neardiag},
\begin{equation}\label{fNexpand}f_N(t)=e^{-|v_N|^2t^2}\left[1+N^{-1/2}\wt
R_{N}(tv_N)\right]\,,\end{equation} where
$$\wt R_{N}(v)=R_1(P_N;0,v,N)+R_1(P_N;v,0,N) - R_1(P_N;v,v,N)
-R_1(P_N;0,0,N) +O(N^{-1/2}) \,.$$
The estimate for $R_1$ yields:
\begin{equation}\label{fNerror} \|\wt R_{N}\|_{\ccal^2\{|v|\le C\}}=
O(1)\,\end{equation}
Since $f_N(1)=1$, it follows from (\ref{fNexpand})--(\ref{fNerror}) that
$|v_N|^2=O(N^{-1/2})$.  (A more careful analysis shows that we can replace
$N^{-1/2}$ with $N^{-1}$ in (\ref{fNexpand}) and thus $|v_N|=O(N^{-1/2})$.)

Write $e^x=1+x+x^2\phi(x)$.  We then have
$$f_N(t)=1-|v_N|^2t^2+|v_N|^4 t^4 \phi(|v_N|^2t^2)+ N^{-1/2}\wt R_N(tv_N)
\left[1-|v_N|^2t^2+|v_N|^4 t^4 \phi(|v_N|^2t^2)\right]\,.$$
Thus by (\ref{fNerror}), $$ f_N''(t)= -2|v_N|^2 +O(|v_N|^4)
+O(N^{-1/2}|v_N|^2)\,,\quad |t|\le 1\,.$$  Since $|v_N|=o(1)$, it follows
that $$0=f_N''(t_N)=(-2+o(1))|v_N|^2\,,$$ which contradicts the assumption
that $v_N\ne 0$.
\end{proof}

\section{Asymptotically  holomorphic versus almost holomorphic
sections}\label{s-ah}

We now use the scaling asymptotics of Theorem \ref{neardiag} to prove Theorem
\ref{aah}, which states that $\nu_{\infty}$-almost every sequence $\{s_N\}$ of
sections (with unit $\lcal^2$-norm) satisfies the  sup-norm estimates
\begin{eqnarray*}\|s_N\|_\infty + \|\bar{\partial} s_N\|_\infty &=& 
O(\sqrt{\log N})\,,\\
\|\nabla^k s_N\|_\infty + \|\nabla^k\bar{\partial} s_N\|_\infty &=&
O(N^{\frac{k}{2}}\sqrt{\log N})\,,\end{eqnarray*} for $k=1,2,3,\dots$.

The following elementary probability lemma is central to our arguments:

\begin{lem}\label{calculus} Let $A\in S^{2d-1}\subset\C^d$, and give
$S^{2d-1}$ Haar probability measure.  Then the probability
that a random point $P\in S^{2d-1}$ satisfies the bound
$|\langle
P,A\rangle|>\la$ is $(1-\la^2)^{d-1}$. \end{lem}

\begin{proof} We can assume without loss of generality that
$A=(1,0,\ldots,0)$.  Let
$$V_\la=\vol \big(\{P\in
S^{2d-1}:|P_1|>\la\}\big)\qquad (0\le\la<1)\,,$$ where $\vol$
denotes $(2d-1)$-dimensional Euclidean volume.  Our desired probability
equals $V_\la/V_0$. Let $\sigma_n=\vol (S^{2n-1})=\frac{2\pi^n}{(n-1)!}$. We
compute
\begin{eqnarray*}V_\la
&=& \int_\la^1\sigma_{d-1}(1-r^2)^{\frac{2d-3}{2}}\frac{2\pi
rdr}{\sqrt{1-r^2}}\ =\ 2\pi\sigma_{d-1}\int_\la^1(1-r^2)^{d-2}rdr\\&=&
\frac{\pi\sigma_{d-1}}{d-1}(1-\la^2)^{d-1}\ =\ \sigma_d
(1-\la^2)^{d-1}\,.\end{eqnarray*}
Therefore  $V_\la/V_0=(1-\la^2)^{d-1}$. \end{proof}

\subsection{Notation}\label{s-notation}
For the readers' convenience, we summarize here our notation for the various
differential operators that we use in this section and elsewhere in the paper:

\begin{enumerate}
\item[a)] Derivatives on $M$:
\begin{itemize}
\item $\frac{\d}{\d z_j} =
\half \frac{\d}{\d x_j}-\frac{i}{2}\frac{\d}{\d y_j}\,,\quad
\frac{\d}{\d\bar
z_j} = \half \frac{\d}{\d x_j}+\frac{i}{2}\frac{\d}{\d y_j}$;\\[-4pt]
\item $Z_j^M=\frac{\d}{\d z_j} +\sum \bar B_{jk}(z)\frac{\d}{\d\bar z_k},\
\bar Z_j^M=\frac{\d}{\d\bar z_j} +\sum  B_{jk}(z)\frac{\d}{\d z_k},\
B_{jk}(P_0)=0$,\\ $\{Z_1,\dots,Z_m\}$ is a local frame for $T^{1,0}M$.
\\[-4pt]\end{itemize}
\item[b)] Derivatives on $X$:
\begin{itemize}
\item $\frac{\d^h}{\d z_j} = \frac{\d}{\d
z_j} -A_j(z)\frac{\d}{\d
\theta} =$ horizontal lift of $\frac{\d}{\d z_j}$, \ $A_j(P_0)=0$;\\[-4pt]
\item $Z_j =$ horizontal lift of $Z_j^M$;\\[-4pt]
\item $d^h=\d_b+ \dbar_b$ = horizontal exterior derivative on $X$.\\[-4pt]
\end{itemize}
\item[c)] Covariant derivatives on $M$:
\begin{itemize}
\item $\nabla:\ccal^\infty(M,L^N\otimes (T^*M)^{\otimes k})\to
\ccal^\infty(M,L^N\otimes (T^*M)^{\otimes (k+1)})$;\\[-4pt]
\item
$\nabla^k=\nabla\circ \cdots\circ \nabla:\ccal^\infty(M,L^N)\to
\ccal^\infty(M,L^N\otimes (T^*M)^{\otimes k} )$;\\[-4pt]
\item $\nabla=\d +\dbar,\ \  \dbar:\ccal^\infty(M,L^N)\to
\ccal^\infty(M,L^N\otimes T^{*0,1}M)$.\\[-4pt] 
\end{itemize}
\item[d)] Derivatives on $X\times X$:
\begin{itemize}
\item $d_j^1, d_j^2$: the operator $\frac{\d^h}{\d z_j}$ applied to the first
and second factors, respectively;\\[-4pt]
\item $Z_j^1, Z_j^2$: the operator $Z_j$ applied to the first
and second factors, respectively.
\end{itemize}
\end{enumerate} 

\subsection{The estimate $\|s_N\|_{\infty}/\|s_N\|_2 = O(\sqrt{\log
N})$ almost surely}

Throughout this section we assume that $\|s_N\|_{\lcal^2} = 1$.
We begin the proof of Theorem
\ref{aah} by showing that
\begin{equation} \nu_N\left\{s^N\in SH^0_J(M,L^N): \sup_M|s_N|>C\sqrt{\log
N}\right\} < O\left(\frac{1}{N^2}\right)\,,
\label{largesupnorm}\end{equation}
for some constant $C<+\infty$.  (In fact, for any $k>0$, we can bound the
probabilities by $O(N^{-k})$ by choosing $C$ to be sufficiently large.)
The
estimate (\ref{largesupnorm}) immediately implies that
$$\limsup_{N\to\infty}\frac{\sup_X|s^N|}{\sqrt{\log N}} \le C\qquad
\mbox{\rm
almost surely}\,,$$ which gives the first statement of Theorem
\ref{neardiag}.

We now show (\ref{largesupnorm}), following an approach inspired by
Nonnenmacher and Voros \cite{NV}. Recalling (\ref{lift}), we note that
\begin{equation}\label{coherentstate} \Pi_N(x,y)=\sum_{j=1}^{d_N}
S^N_j(x)\overline{S^N_j(y)}=\langle
\tilde\Phi_N(x),\tilde\Phi_N(y)\rangle\,.
\end{equation}
Let $s^N=\sum_{j=1}^{d_N}c_jS^N_j\ $ ($\sum|c_j|^2=1$) denote a random
element of $SH^0(M,L^N)=S\hcal^2_N(X)$, and write $c=(c_1,\dots,c_{d_N})$.
Recall that
\begin{equation}\label{sNx}s^N(x)=\int_X\Pi_N(x,y)s^N(y)dy=
\sum_{j=1}^{d_N}c_j  S^N_j(x)=
c\cdot\tilde\Phi_N(x)\,.\end{equation} Thus
\begin{equation}\label{cos}|s^N(x)|=\|\tilde\Phi_N(x)\|\cos \theta_x\,,\quad
\mbox{where\ } \cos \theta_x = \frac{\left| c\cdot
\tilde\Phi_N(x)\right|}{\| \tilde\Phi_N(x)\|} \,.\end{equation} (Note that
$\theta_x$ can be interpreted as the distance in $\CP^{d_N-1}$ between
$[\bar c]$ and
$\tilde\Phi_N(x)$.)

We have by Theorem
\ref{tyz}(a),
\begin{equation}\label{tyza}\|\tilde\Phi_N(x)\|=\Pi_N(x,x)^\half
= N^{m/2}
+O(N^{m/2-1})=(1+\ep_N)N^{m/2} \,,\end{equation} where $\ep_N$ denotes a
term
satisfying the uniform estimate
\begin{equation}\label{epN} \sup_{x\in X}|\ep_N(x)| \le
O\left(\frac{1}{N}\right)\,.\end{equation}
Now fix a point $x\in X$.  By Lemma
\ref{calculus}, \begin{eqnarray}\label{prob1}
\nu_N\left\{s^N:\cos\theta_x\ge C N^{-m/2}\sqrt{\log
N}\right\}&=& \left(1-\frac{C^2\log N}{N^m}\right)^{d_N-1}
\nonumber \\ &\le & \left(e^{-\frac{C^2\log
N}{N^m}}\right)^{d_N-1}\ =\ N^{-C^2 N^{-m}(d_N-1)}\;.\end{eqnarray}

We can cover $M$ by a
collection of $k_N$ balls $B(z^j)$ of radius
\begin{equation}\label{RN} R_N:=\frac{1}{N^{\frac{m+1}{2}}}\end{equation}
 centered at points
$z^1,\dots,z^{k_N}$, where
$$k_N \le O(R^{-2m})\le
O(N^{m(m+1)})\,.$$
By  (\ref{prob1}), we have
\begin{equation} \label{centers} \nu_N \left\{s^N\in
SH^0_J(M,L^N):\max_{j} \cos\theta_{x^j} \ge C N^{-m/2}\sqrt{\log
N}
\right\}\le k_N N^{-C^2 N^{-m}(d_N-1)}\;,\end{equation}
where $x^j$ denotes a point in $X$ lying above $z^j$.

Equation (\ref{centers}) together with (\ref{cos})--(\ref{tyza})
implies (by the argument below) that the
desired sup-norm estimate holds at the centers of the small balls with high
probability. To complete the proof of (\ref{largesupnorm}), we first need
to extend (\ref{centers}) to points within the balls.
To do this,  we consider an arbitrary point $w^j\in B(z^j)$, and choose
points
$y^j\in X$ lying above  the points $w^j$.  We must estimate the
distance, which we denote by $\de_{N}^j$, between $\Phi_N(z^j)$ and
$\Phi_N(w^j)$ in
$\CP^{d_N-1}$. Letting $\ga$ denote the geodesic in
$M$ from $z^j$ to $w^j$, we conclude by Theorem \ref{tyz}(b) that
\begin{eqnarray}\de_{N}^j &\le & \int_{\Phi_{N*}\ga}\sqrt{\om_{FS}}
\ =\ \int_\ga \sqrt{\Phi_{N}^*\om_{FS}}
\ \le \ \sqrt{N}\int_\ga
(1+\ep_N)\big)\sqrt{\om}\nonumber \\
& \le& (1+\ep_N)N^\half  R_N
\ =\  \frac{1+\ep_N}{N^{m/2}}\,.
\label{delta} \end{eqnarray}

By the triangle inequality
in $\CP^{d_N-1}$, we have
$|\theta_{x^j}-\theta_{y^j}|\le\de_N^j$. Therefore by (\ref{delta}),
\begin{equation}\label{triangle}\cos\theta_{x^j} \ge
\cos\theta_{y^j}-\de_N^j
\ge
\cos\theta_{y^j}- \frac{1+\ep_N}{N^{m/2}}\,.\end{equation}
By (\ref{triangle}),
$$\cos\theta_{y^j} \ge \frac{(C+1)\sqrt{\log N}}{N^{m/2}}\Rightarrow
\cos\theta_{x^j} \ge \frac{(C+1)\sqrt{\log N} -(1+\ep_N)}{N^{m/2}}
\ge  \frac{C\sqrt{\log N}}{N^{m/2}}$$ and thus
$$\begin{array}{l}\left\{s^N\in
SH^0_J(M,L^N):\sup \cos\theta \ge (C+1) N^{-m/2}\sqrt{\log
N}
\right\}\\[8pt] \qquad\qquad\qquad \subset \left\{s^N\in
SH^0_J(M,L^N):\max_{j} \cos\theta_{x^j} \ge C N^{-m/2}\sqrt{\log
N}
\right\}\,.\end{array}$$

Hence
by  (\ref{centers}),
\begin{equation} \label{cos*} \nu_N \left\{s^N\in
SH^0_J(M,L^N):\sup \cos\theta \ge (C+1) N^{-m/2}\sqrt{\log
N}
\right\}\le k_N N^{-C^2 N^{-m}(d_N-1)}\;.\end{equation}

By the
Riemann-Roch formula of Boutet de Monvel - Guillemin \cite[\S 14]{BG} (which
is a consequence of Theorem \ref{COMPLEX}), we have the estimate for the
dimensions $d_N$:
\begin{equation}\label{RR}
d_N=\frac{c_1(L)^m}{m!}N^m+O(N^{m-1})\,.\end{equation}
We can also obtain (\ref{RR}) from Theorem
\ref{tyz}(a) as follows: We note first that
$$\int_X \Pi_N(x,x)d\vol_X=\int_X \sum_{j=1}^{d_N}|S_j^N(x)|^2d\vol_X
=d_N\,.$$ On the other hand, by Theorem \ref{tyz}(a),
$$\int_X \Pi_N(x,x)d\vol_X =
\left[\frac{1}{\pi^m}N^m+O(N^{m-1})\right]\vol(X)\,,$$ where $$\vol(X)
=\vol(M)=\int_M\frac{1}{m!}\om^m = \frac{\pi^m}{m!}c_1(L)^m\,.$$
Equating the above computations of the integral yields (\ref{RR}).

It follows from (\ref{cos}), (\ref{tyza}), (\ref{cos*}) and (\ref{RR}) that
\begin{eqnarray*}\nu_N \left\{s^N\in SH^0_J(M,L^N):\sup_M |s^N|\ge  (C+2)
\sqrt{\log N}
\right\}\hspace{1.25in}\\
\le k_N N^{-C^2 N^{-m}(d_N-1)}\le O\left(N^{m(m+1)-\frac{C^2}
{m!+1}}\right)\;.\end{eqnarray*}
Choosing $C=(m+1)\sqrt{m!+1}$, we obtain
(\ref{largesupnorm}).

\begin{rem} An alternate proof of this estimate, which does not depend on
Tian's theorem, is given by the case $k=0$ of the $\ccal^k$ estimate in the
next section.
\end{rem}

\subsection{The estimate $\|\nabla^k s_N\|_{\infty}/\|s_N\|_2 =
O(\sqrt{N^k \log N})$ almost surely} \label{almostsurelyCkbounded}

The proof of this and the other assertions of Theorem \ref{aah} follow the
pattern of the
above sup-norm estimate.
First we note a consequence (Lemma \ref{horiz-bound2}) of our near-diagonal
asymptotics.  Recall that a differential operator on $X$ is horizontal if
it is generated by horizontal vector fields.  In particular the operators
$\nabla^k:\ccal^\infty(M,L^N)\to \ccal^\infty(M,L^N\otimes (T^*M)^{\otimes
k})$ are given by (vector valued) horizontal differential operators
(independent of $N$) on
$X$. By definition, horizontal differential operators on
$X\times X$ are generated by the horizontal differential operators on the
first and second factors. We begin with the following estimate:

\begin{lem}\label{horiz-bound1} Let $P_k$ be a horizontal differential
operator of order $k$ on $X\times X$. Then
$$P_k\Pi_N(x,y)|_{x=y}=O(N^{m+k/2})\,.$$\end{lem}

\begin{proof}
Let $x_0=(P_0,0)$ be an arbitrary point of $X$, and choose local {\bf
real} coordinates $(x_1,\dots,x_{2m},\theta)$ about $(P_0,0)$ as in the
hypothesis of Theorem \ref{neardiag} (with $z_q=x_q+ix_{m+q}$).
We let $\frac{\d^h}{\d x_q}$ denote the
horizontal lift of $\frac{\d}{\d x_q}$ to $X$:
$$\frac{\d^h}{\d x_q} = \frac{\d}{\d x_q} -\tilde A_q(x)\frac{\d}{\d
\theta}\,,\quad \tilde 
A_q =(\al, \frac{\d}{\d x_q})\,.$$ Since $\left.\frac{\d}{\d
x_q}\right|_{x_0}$ is assumed to be horizontal, we have $\tilde A_q(P_0)=0$.

We let $d_q^1, d_q^2$ denote the operator $\frac{\d^h}{\d x_q}$ applied to
the first and second factors, respectively, on $X\times X$.
For this result, we need only the zeroth order estimate of Theorem
\ref{neardiag}:
\begin{equation}\label{zeroth}
\Pi_N(\frac{u}{\sqrtn},\frac{s}{N};
\frac{v}{\sqrtn},\frac{t}{N}) = {N^m}e^{i(s-t)+\psi_2(u,v)}
\rcal(P_0,u,v,N)\,,\end{equation}
where $\rcal(P_0,u,v,N)$ denotes a term satisfying
the remainder estimate of Theorem
\ref{neardiag}: $\|\rcal(P_0,u,v,N)\|_{\ccal^j(\{|u|\le
\rho,\ |v|\le
\rho\}}\le C_{j,\rho}$ for $j\ge 0,\,\rho>0$, where $C_{j,\rho}$ is
independent of the point
$P_0$ and choice of coordinates.

Differentiating (\ref{zeroth})  and noting that $\d/\d x_q = \sqrtn
\d/\d u_q$, $\d/\d\theta = N\d/\d s$, we have
\begin{eqnarray}\lefteqn{d_q^1\Pi_N(\frac{u}{\sqrtn},\frac{s}{N};
\frac{v}{\sqrtn},\frac{t}{N})} \nonumber \\
&=& \sqrtn \left( \frac{\d}{\d
u_q}-\sqrtn\tilde A_q(P_0+{\textstyle{\frac{u}{\sqrtn}}})\frac{\d}{\d
s}\right)
\left( {N^m}e^{i(s-t)+\psi_2(u,v)}\rcal\right) \nonumber\\
&=&{N^{m+1/2}}
e^{i(s-t)+\psi_2(u,v)}\left\{\left[L_q(u,v)-i\sqrtn \tilde
A_q({\textstyle{\frac{u}{\sqrtn}}})\right] +\frac{\d}{\d u_q}\rcal
\right\}
\nonumber\\ &=& {N^{m+1/2}}
e^{i(s-t)+\psi_2(u,v)} \wt\rcal \ = \
O(N^{m+1/2})\,,\label{expand0}
\end{eqnarray} where $L_q:=\frac{\d\psi_2}{\d u_q}$ is a linear function.
The same estimate holds for $d^2_q\Pi_N$.  Indeed, the above computation
yields:
\begin{equation}\label{estimate0} d_q^je^{i(s-t)+\psi_2(u,v)}
\rcal(P_0,u,v,N) = \sqrtn
e^{i(s-t)+\psi_2(u,v)} \wt \rcal(P_0,u,v,N)\,,\end{equation} for $j=1,2,\
q=1,\dots,2m$.  The desired estimate follows by iterating
(\ref{estimate0}). \end{proof}

\begin{rem} The assumption that $P_k$ is horizontal in Lemma
\ref{horiz-bound1} is necessary, since the operator $\d/\d\theta$
multiplies the estimate by $N$ instead of $\sqrtn$. \end{rem}

\begin{lem} Let $P_k$ be a horizontal differential operator of
order
$k$ on
$X$. Then
$$ \sup_X\|P_k \tilde\Phi_N\| = O(N^\frac{m+k}{2})\,.$$
\label{horiz-bound2}\end{lem}

\begin{proof} Let $P_k^1,\ P_k^2$  denote the operator $P_k$ applied to
the first and second factors, respectively, on $X\times X$. Differentiating
(\ref{coherentstate}) and restricting to the diagonal, we obtain
\begin{equation}\label{coherentstate1}
P_k^1 \bar P_k^2\Pi_N(x,x)= \left
\|P_k\tilde\Phi_N(x)\right\|^2\,.
\end{equation}  The conclusion follows from (\ref{coherentstate1}) and Lemma
\ref{horiz-bound1} applied to the horizontal differential operator (of order
$2k$)
$P_k^1
\bar P_k^2$ on $X\times X$.
\end{proof}

We are now ready to use the small-ball method of the previous section to
show that $\|\nabla^k s_N\|_{\infty}/\|s_N\|_{\lcal^2} = O(\sqrt{N^k \log
N})$ {\it almost surely.\/}  It is sufficient to show that
\begin{equation} \nu_N\left\{s^N\in SH^0_J(M,L^N): \sup_M|\nabla^k
s_N|>C\sqrt{N^k   \log N}\right\} <
O\left(\frac{1}{N^2}\right)\,,
\label{largeCknorm}\end{equation} for $C$ sufficiently large.
To verify (\ref{largeCknorm}), we may regard $s_N$ as a function on
$X$ and replace
$\nabla^k$ by a horizontal $r_\theta$-invariant differential operator of
order
$k$ on
$X$.

As before, we let $s^N=\sum c_j s^N_j$ denote a random element of
$S\hcal^2_N(X)$.
By (\ref{sNx}), we have
\begin{equation}\label{dqsn}P_ks^N(x)=\int_X
P_k^1\Pi_N(x,y)s^N(y)dy= \sum_{j=1}^{d_N}c_j P_kS^N_j(x) = c\cdot
P_k\tilde\Phi_N(x)
\,.\end{equation} We then have
\begin{equation}\label{cos1}|P_ks^N(x)|=\|P_k\tilde\Phi_N(x)\|\cos
\theta_x\,,\quad
\mbox{where\
} \cos \theta_x = \frac{\left| c\cdot
P_k\tilde\Phi_N(x)\right|}{\|
P_k\tilde\Phi_N(x)\|} \,.\end{equation}
Now fix a point $x\in X$. As before, (\ref{prob1}) holds, and hence by Lemma
\ref{horiz-bound2} we have
\begin{equation} \label{supPk} \nu_N \left\{s^N\in S\hcal^2_N:
|P_ks^N(x)| \ge C'
\sqrt{N^k\log N}
\right\}\le k_NN^{-C^2 N^{-m}(d_N-1)}\;,\end{equation}
where $C'=C \sup_{N,x}N^{-(m+k)/2}|P_k\tilde\Phi_N(x)|$.

We again  cover $M$ by a collection of $k_N$ very small balls $B(z^j)$ of
radius
$R_N=N^{-\frac{m+1}{2}}$ and first show that the probability of the
required condition holding at the centers of all the balls is small.
Choosing points
$x^j\in X$ lying above the centers $z^j$ of the balls, we
then have
\begin{equation} \label{centers2} \nu_N \left\{s^N\in S\hcal^2_N:\max_{j}
|P_ks^N(x_j)| \ge C'
\sqrt{N^k\log N}
\right\}\le k_NN^{-C^2 N^{-m}(d_N-1)}\;.\end{equation}
Now suppose that $w^j$ is an arbitrary point in $B(z^j)$, and let $y^j$ be
the point of $X$ above $w^j$ such that the horizontal lift of the geodesic
from
$z^j$ to $w^j$ connects $x^j$ and $y^j$.  Hence by Lemma
\ref{horiz-bound2}, we have
\begin{equation}\label{dschange}\|P_k\tilde
\Phi_N(x^j)-P_k\tilde
\Phi_N(y^j)\|\le \sup_M\|d^h(P_k\tilde
\Phi_N)\|
r_N = O(N^{\frac{m+k+1}{2}})r_N=O(N^{\frac{k}{2}})\,.\end{equation}
It follows as before from (\ref{centers2}) and
(\ref{dschange}) that
\begin{eqnarray*} \nu_N \left\{s^N\in S\hcal^2_N:\sup_X
|P_ks^N| \ge (C'+1)
\sqrt{N^k\log N}
\right\}\hspace{1.25in}\\
\le k_N N^{-C^2 N^{-m}(d_N-1)} \le
O\left(N^{m(m+1)-\frac{C^2}
{m!+1}}\right)\;.\end{eqnarray*}
(Here, we used the fact that $|P_ks^N|$ is constant on the fibers of
$\pi:X\to M$.) Thus, (\ref{largeCknorm})  holds with $C$ sufficiently
large.\qed

\subsection{The estimate $\|\bar{\partial}
s_N\|_{\infty}/\|s_N\|_2 =
O(\sqrt{\log N})$ almost surely}

The proof of the $\dbar s_N$ estimate follows the
pattern of the
above estimate.  However, there is one crucial difference: we must
show the following upper bound for the modulus of $\dbar_b\tilde\Phi_N$.
This estimate is a factor of $\sqrtn$ better than the one for $d^h
\tilde\Phi_N$ arising from Lemma
\ref{horiz-bound2}; the proof depends on the precise second order
approximation of Theorem
\ref{neardiag}.

\begin{lem}\label{tyzabar} $\qquad\sup_X\|\dbar_b\tilde\Phi_N(x)
\|\le O(N^{m/2})
\,.$\end{lem}
\begin{proof}
Let $x_0=(P_0,0)$ be an arbitrary point of $X$, and choose preferred local
coordinates $(z_1,\dots,z_{m},\theta)$ about $(P_0,0)$ as in the
hypothesis of Theorem \ref{neardiag}.  We lift a local frame $\{\bar
Z_q^M\}$ of the form (\ref{localtangentframe}) to obtain the local frame
$\{\bar Z_1,\dots,\bar Z_m\}$ for
$H^{0,1}X$ given by
\begin{equation}\label{barz}\bar Z_q = \frac{\d^h}{\d\bar z_q} +
\sum_{r=1}^m B_{qr}(z)\frac{\d^h}{\d z_r}\,,\qquad
B_{qr}(P_0)=0\,.\end{equation}
It suffices to show that
\begin{equation}\label{barzh}N^{-m/2}|\bar Z_q \tilde \Phi_N (x_0)| \le
C\,,\end{equation} where $C$ is a constant independent of $x_0$.

By Theorem \ref{neardiag}, we have \begin{eqnarray}
\lefteqn{N^{-m} \Pi_N(\frac{u}{\sqrtn},\frac{s}{N};
\frac{v}{\sqrtn},\frac{t}{N})}\nonumber \\& =&
\frac{1}{\pi^m}\phi_0(u,s)\overline{\phi_0(v,t)}e^{u\cdot
\bar v}
\left[1 +
\frac{1}{\sqrt{N}}
b_1(P_0,u,v)+\frac{1}{N}R_2(P_0,u,v,N)\right]
\,,\label{asymptotics1}
\end{eqnarray} where 
$$\phi_0(z,\theta)=e^{i\theta-|z|^2/2}\,.$$
(The function $\phi_0$ is the `ground state' for the `annihilation
operators'
$\bar Z_q$ in the Heisenberg model; see the remark in \S
\ref{s-heisenberg} and
\cite[\S 1.3.2]{BSZ2}).  In our case, $\bar Z_q\phi_0$ does not
vanish as in the model case, but instead  satisfies the
asymptotic bound (\ref{groundstate}) below.)   

Recalling (\ref{dhdzj}) and (\ref{dhoriz}), we have
\begin{equation}\frac{\d^h}{\d\bar z_q} = \frac{\d}{\d\bar z_q}
+\left[ -\frac{i}{2}z_q-R_1^{\bar A_q}(z)\right]\frac{\d}{\d \theta}
\,,\label{dh}\end{equation} where
$R_1^{\bar A_q}(z)=O(|z|^2)$.  Recalling that $z=u/\sqrtn$,
$\theta=s/N$, we note that $\phi_0(u,s)=e^{iN\theta
-N|z|^2/2}=\phi_0(z,\theta)^N$, 
and thus by (\ref{dh}),
\begin{equation}\label{groundstate} \frac{\d^h}{\d\bar z_q}
\phi_0 (u,s)=\frac{\d^h}{\d\bar z_q}e^{iN\theta
-N|z|^2/2}  = -iN
R_1^{\bar A_q}(\frac{u}{\sqrtn})
\phi_0(u,s) = \rcal(P_0,u,N)\phi_0(u,s)\,,\end{equation}
where as before $\rcal$ denotes a term satisfying the remainder estimate
of Theorem \ref{neardiag}.

We let
$ Z_q^1,  Z_q^2$ denote the operator
$ Z_q$ applied to the first and second factors, respectively,
on $X\times X$; we similarly let $d_q^1,
d_q^2$ denote the operator $\frac{\d^h}{\d z_q}$ applied to the factors of
$X\times X$. Equation (\ref{coherentstate1}) tells us that
\begin{equation}\label{barzh2}
\|\bar Z_q\tilde\Phi_N(x)\|^2=
\bar Z_q^1 Z_q^2 \Pi_N(x,x)\,.
\end{equation}
By (\ref{barz}), \begin{equation}\label{4terms}\bar Z_q^1 Z_q^2
=  \left(\overline{d^1_q}+\sum_{r=1}^mB_{qr}(z)
d^1_r\right)
\left(d^2_q+\sum_{\rho=1}^m \bar B_{q\rho}(w)
\overline{d^2_\rho}\right)
\,,\end{equation} where we recall that $B_{qr}(P_0)=0$.

Differentiating
(\ref{asymptotics1}), again noting that
$\d/\d z_q =\sqrtn\d/\d u_q$, $\d/\d w_q =\sqrtn\d/\d v_q$ and using
(\ref{groundstate}), we obtain
\begin{eqnarray} \ \lefteqn{N^{-m}\left(\overline{ d^1_q} d^2_q
\Pi_N \right)(\frac{u}{\sqrtn},\frac{s}{N};
\frac{v}{\sqrtn},\frac{t}{N})}\nonumber \\
&\quad =&\frac{1}{\pi^m}
\phi_0(u,s)\overline{\phi_0(v,t)}e^{u\cdot
\bar v} \left[\sqrtn \frac{\d^2}{\d\bar u_q\d v_q}b_1 +\wt\rcal\right]
\,.\label{asymptotics2}\end{eqnarray}
Since $b_1$ has no terms that are quadratic in $(u,\bar u,v,\bar v)$, it
follows from (\ref{4terms})--(\ref{asymptotics2}) that
\begin{equation} N^{-m}\left|\bar Z_q^1 Z_q^2
\Pi_N(P_0,0;P_0,0)\right| =
N^{-m}\left|\overline{ d^1_q} d^2_q
\Pi_N(P_0,0;P_0,0)\right| = \frac{1}{\pi^m}\left|
\wt\rcal(P_0,0,0,N)\right| \le O(1)\,. \label{barzh1}\end{equation}
The desired estimate (\ref{barzh}) now follows immediately from
(\ref{barzh2}) and (\ref{barzh1}). \end{proof}

By covering $M$ with small
balls and repeating the argument of the previous section, using Lemma
\ref{tyzabar}, we conclude that
\begin{equation}
\nu_N\left\{s^N\in SH^0_J(M,L^N):
\sup_M|\dbar s_N|>C\sqrt{\log N}\right\} < O\left(\frac{1}{N^2}\right)\,.
\end{equation}
Thus $\|\bar{\partial}
s_N\|_{\infty}/\|s_N\|_2 =
O(\sqrt{\log N})$ almost surely.

\subsection {The estimate $\|\nabla^k\bar{\partial}
s_N\|_{\infty}/\|s_N\|_2 =
O(\sqrt{N^k\log N})$ almost surely}

To obtain this final estimate of Theorem \ref{aah}, it
suffices to verify the probability
estimate
\begin{equation} \nu_N\left\{s^N\in SH^0_J(M,L^N): \sup_M|\nabla^k \dbar
s_N|>C\sqrt{N^k   \log N}\right\} <
O\left(\frac{1}{N^2}\right)\,.
\label{largeCknorm2}\end{equation}
Equation (\ref{largeCknorm2}) follows by again repeating the argument of
\S \ref{almostsurelyCkbounded}, using the following lemma.

\begin{lem} Let $P_k$ be a horizontal differential operator of
order
$k$ on
$X$ ($k\ge 0$). Then
$$ \sup_X|P_k \dbar_b\tilde\Phi_N| = O(N^\frac{m+k}{2})\,.$$
\label{horiz-bound3}\end{lem}

\begin{proof} It suffices to show that
\begin{equation}\label{supkq}\sup_U|P_k \bar Z_q^k\tilde\Phi_N| =
O(N^\frac{m+k}{2})\end{equation} for a local frame $\{\bar Z_q\}$ of
$T^{0,1}M$ over
$U$. As before, we have
\begin{equation}\label{coherentstate2}
P_k^1  \bar P_k^2 \bar Z_q^1 Z_q^1 \Pi_N(x,x)= \left
|P_k\bar Z_q^{h}\tilde\Phi_N(x)\right|^2\,.
\end{equation}

We claim that
\begin{equation}\label{step0} N^{-m}\bar Z^1_q Z^2_q
 \Pi_N = \frac{1}{\pi^m}
e^{i(s-t)+\psi_2(u,v)}\left[\sqrtn  \frac{\d^2}{\d\bar u_q\d v_q}b_1
+\rcal(P_0,u,v,N)\right]\,.\end{equation}
To obtain the estimate (\ref{step0}), we recall from (\ref{4terms}) in
the proof of Lemma \ref{tyzabar} that
\begin{equation}\label{4terms*}\bar Z_q^1 Z_q^2 =
\overline{d^1_q}d^2_q+
\sum_{\rho=1}^m \bar B_{q\rho}(w)\overline{d^1_q}
\overline{d^2_\rho} + \sum_{r=1}^mB_{qr}(z) d^1_r d^2_q +
\sum_{r,\rho}B_{qr}(z)\bar B_{q\rho}(w)  d^1_r\overline{d^2_\rho}\,.
\end{equation}
Equation (\ref{asymptotics2}) says that the first term of
$N^{-m}\bar Z^1_q Z^2_q
\Pi_N$
coming from the expansion (\ref{4terms*}) satisfies the estimate of
(\ref{step0}).  To obtain the estimate for the second term, we compute:
\begin{eqnarray}
\lefteqn{N^{-m}\overline{d^2_\rho}
\Pi_N(\frac{u}{\sqrtn},\frac{s}{N};
\frac{v}{\sqrtn},\frac{t}{N})\ =\ \frac{\sqrtn}{\pi^m}
e^{i(s-t)+\psi_2(u,v)}}\nonumber \\
& & \cdot  \left(\left[\frac{\d\psi_2}{\d \bar v_\rho}+i\sqrtn
A_\rho(\frac{v}{\sqrtn})\right] \left[1+ \frac{1}{\sqrtn}b_1+
\frac{1}{N}R_2
\right] +\frac{1}{\sqrtn}\frac{\d\ b_1}{\d\bar v_\rho}
+\frac{1}{N} \frac{\d\ R_2}{\d\bar v_\rho}\right)\nonumber \\
&=& \frac{\sqrtn}{\pi^m}
\phi_0(u,s)\overline{\phi_0(v,t)}e^{u\cdot
\bar v} \left[\frac{\d\psi_2}{\d \bar
v_\rho} +L_\rho(v)+\frac{1}{\sqrtn}\wt\rcal\right]\,,\end{eqnarray}
where $ L_{\rho}$ is a linear function. Since $\d^2\psi_2/\d\bar
u_q \d\bar v_\rho \equiv 0$, it then follows that
\begin{equation} N^{-m}\overline{ d^1_q}\overline{d^2_\rho}
 \Pi_N = \frac{\sqrtn}{\pi^m}
e^{i(s-t)+\psi_2(u,v)} \frac{\d}{\d \bar
u_q}\wt\rcal(P_0,u,v,N)\,.\label{mixedterm}\end{equation}
The estimate (\ref{step0}) for the second term follows from
(\ref{mixedterm}), using the fact that
$B_{q\rho}(\frac{v}{\sqrtn})=\frac{1}{\sqrtn}L_{q\rho}(v)+\cdots$.
The proofs of the estimate for the third and fourth terms are
similar.

The desired estimate (\ref{supkq}) follows as before from
(\ref{coherentstate2}), (\ref{step0}), and (\ref{estimate0}), using the
fact that $ \frac{\d^2}{\d\bar
u_q\d v_q}b_1$ is linear.
\end{proof}

\section {The joint probability distribution}\label{s-universality}
In this section, we shall use Theorem \ref{neardiag} and the methods of
of \cite{BSZ2} to prove Theorem \ref{usljpd-sphere} and its analogue for
Gaussian measures (Theorem \ref{usljpd}), which say that the joint
probability distributions on almost complex symplectic manifolds have the
same universal scaling limit as in the complex case.

\subsection{Generalized Gaussians}\label{s-gaussians}

Recall that a Gaussian measure on $\R^n$ is a measure of the form
$$\ga_\De = \frac{e^{-\half\langle\De\inv
x,x\rangle}}{(2\pi)^{n/2}\sqrt{\det\De}}  dx_1\cdots dx_n\,,$$ where $\De$
is a positive definite symmetric $n\times n$ matrix.  The matrix $\De$ gives
the second moments of $\ga_\De$:
\begin{equation}\label{moments}\langle x_jx_k
\rangle_{\ga_\De}=\De_{jk}\,.\end{equation} This
Gaussian measure is
also characterized by its Fourier transform
\begin{equation}\label{muhat}\wh{\ga_\De}(t_1,\dots,t_n) = e^{-\half\sum
\De_{jk}t_jt_k}\,.\end{equation}
If we let $\De$ be the $n\times n$ identity matrix, we obtain the standard
Gaussian measure on $\R^n$,
$$\ga_n:=\frac{1}{(2\pi)^{n/2}}e^{-\half |x|^2} dx_1\cdots dx_n\,,$$
with the
property that the
$x_j$ are independent Gaussian variables with mean 0 and variance 1. Hence
$$\langle \|x\|^2\rangle_{\ga_n}= \sum_{j=1}^n\langle
x_j^2\rangle_{\ga_n}=n\,.$$  Since we wish to put Gaussian measures on the
spaces $H^0_J(M,L^N)$ with rapidly growing dimensions, it is useful to
consider the {\it normalized standard Gaussians}
$$\tilde\ga_n:= k_n e^{-\frac{n}{2} |x|^2} dx_1\cdots dx_n\,,\quad k_n=
\left(\frac{n}{2\pi}\right)^{n/2}\,,$$
which have the property that
$$\langle \|x\|^2\rangle_{\tilde\ga_n} =1\,.$$

The push-forward of a Gaussian measure by a surjective linear map is also
Gaussian. In the next section, we shall push forward Gaussian
measures (on the spaces
$H^0_J(M,L^N)$) by linear maps that are sometimes not surjective.  Since
these non-surjective push-forwards are  singular measures, we need to
consider the case where
$\De$ is positive semi-definite.  In this case, we use (\ref{muhat}) to
define a measure $\ga_\De$, which we call a {\it generalized Gaussian\/}. 
If $\De$ has null eigenvalues, then the generalized Gaussian
$\ga_\De$ is a Gaussian measure on the subspace $\Lambda_+\subset\R^n$
spanned by the positive eigenvectors.  (Precisely, $\ga_\De=\iota_*
\ga_{\De|\Lambda_+}$, where $\iota:\Lambda_+\hookrightarrow \R^n$ is the
inclusion.  For the completely degenerate case
$\De=0$, we have
$\ga_\De=\de_0$.)  Of course, (\ref{moments}) holds for semi-positive $\De$.
One useful property of generalized Gaussians is that the push-forward
by a (not necessarily surjective) linear map $T:\R^n\to\R^m$ of a
generalized Gaussian $\ga_\De$ on $\R^n$ 
is a generalized Gaussian on $\R^m$:
\begin{equation}\label{pushgaussian} T_*\ga_\De=\ga_{T\De T^*}\end{equation}
Another useful property of  generalized Gaussians is the following fact:

\begin{lem}\label{continuity}  The map $\De\mapsto\ga_\De$ is a continuous
map from the positive semi-definite matrices to the space of positive
measures on
$\R^n$ (with the weak topology).\end{lem}

\begin{proof} Suppose that $\De^N\to\De^0$. We must show that
$(\De^N,\phi)\to(\De^0,\phi)$ for a compactly supported test function
$\phi$.  We can assume that $\phi$ is $\ccal^\infty$. It then
follows from (\ref{muhat}) that
$$(\ga_{\De^N},\phi)=(\wh{\ga_{\De^N}},\wh\phi)\to
(\wh{\ga_{\De^0}},\wh\phi)=(\ga_{\De^0},\phi)\,.$$\end{proof}

We shall use the following general result relating spherical measures to
Gaussian measures in order to prove Theorem \ref{usljpd-sphere} on
asymptotics of the joint probability distributions for $SH^0_J(M,L^N)$.

\begin{lem}\label{spherical-vs-gaussian}
 Let $T_N: \R^{d_N} \to R^k$, $N=1,2,\dots$, be a sequence of linear maps,
where $d_N
\to \infty$.  Suppose that $\frac{1}{d_N}T_N T_N^*
\to
\Delta$. Then
$T_{N*} \nu_{d_N} \to \gamma_{\Delta}$. \end{lem}

\begin{proof} 
Let $V_N$ be a $k$-dimensional subspace of $\R^{d_N}$ such that
$V_N^\perp\subset \ker T^N$, and let
$p_N:\R^{d_N}\to V_N$ denote the orthogonal projection.  We decompose
$T_N=B_N\circ A_N$, where
$A_N=d_N^{1/2}p_N:\R^{d_N}\to V_N$, and $B_N =
d_N^{-1/2}T_N|_{V_N}:V_N\buildrel{\approx }\over\to\R^k$.
  Write $$A_{N*}\nu_{d_N} =
\al_N\,,\qquad T_{N*}
\nu_{d_N} =
B_{N*}\al_N = \be_N\,.$$  We easily see that (abbreviating $d=d_N$)

\begin{equation}\label{archimedes}\al_N=A_{N*} \nu_d = \psi_d dx\,,\qquad
\psi_d =
\left\{\begin{array}{cl}\frac{\sigma_{d-k}}{\sigma_d d^k} [1 -  \frac{1}{d}
|x|^2]^{(d-k-2)/2} 
\ &{\rm for}\ |x|<\sqrt{d}\\[6pt]
\ 0 & {\rm otherwise}\end{array}\right.,\end{equation}
where $dx$ denotes Lebesgue measure on $V_N$, and
$\sigma_n=\vol(S^{n-1})=
\frac{2\pi^{n/2}}{\Ga(n/2)}$. (The case $k=1, \
d=3$ of (\ref{archimedes}) is  Archimedes' formula \cite{Arc}.) Since
$[1 -   |x|^2/d]^{(d-k-2)/2}
\to e^{-|x|^2/2}$ uniformly on compacta and $\frac{\sigma_{d-k}}{\sigma_d
d^k}\to \frac{1}{(2\pi)^{k/2}}$, we conclude that $\al_N\to\ga_k$. (This is
the Poincar\'e-Borel Theorem; see Corollary \ref{PB} below.) Furthermore,
$$\left(1 -  \frac{1}{d} |x|^2\right)^{(d-k-2)/2}
\le \exp\left(-\frac{d-k-2}{2d}|x|^2\right)\le 
e^{\frac{k+2}{2}}e^{-\half|x|^2}\qquad{\rm
for}\ d\ge k+2\,,\ 
|x|\le\sqrt{d}\,,$$
 and hence \begin{equation}\label{ldc}\psi_{d_N}(x)\le
C_k e^{-|x|^2/2}\,.\end{equation}

Now let $\phi$ be a compactly supported continuous test function on
$\R^k$. We must show that \begin{equation}\label{true}\int \phi
d\be_N\to\int
\phi d\ga_{\De}\,.\end{equation}  Suppose on the contrary that
(\ref{true}) does not hold. After passing to a subsequence, we may assume that 
$\int \phi d\be_N\to c\ne \int
\phi d\ga_{\De}$.  Since the eigenvalues of $B_N$ are bounded, we can
assume  (again taking a subsequence) that $B_N\to B_0$, where
$$B_0 B_0^*=\lim _{N\to\infty}B_NB_N^*=\lim
_{N\to\infty}\frac{1}{d_N} T_N T_N^*=\De\,.$$  Hence,
$$\int_{\R^k} \phi
d\be_N =\int_{V_N} \phi(B_Nx)\psi_{d_N}(x)dx \to \int_{V_N} \phi(B_0
x)\frac{e^{-|x|^2/2}}{(2\pi)^{k/2}} dx=
 \int_{V_N} \phi(B_0 x)d\ga_k(x)\,,$$
where the limit holds by dominated convergence, using (\ref{ldc}).
By (\ref{pushgaussian}), we have $B_{0*}\ga_k=\ga_{B_0B_0^*}=\ga_\De$, and
hence 
$$ \int_{V_N} \phi(B_0 x)d\ga_k(x) = \int_{\R^k}\phi d\ga_{\De}\,.$$
Thus (\ref{true}) holds for the subsequence, giving a contradiction.
\end{proof}

\medskip We note that the above proof
began by establishing the {\it Poincar\'e-Borel Theorem\/} (which is
a special case of the  of Lemma \ref{spherical-vs-gaussian}):

\begin{cor} {\rm  (Poincar\'e-Borel)} Let $P_d:\R^d\to \R^k$ be given by
$P_d(x)=\sqrt{d}(x_1,\dots,x_k)$.  Then $$P_{d*}\nu_d \to \ga_k\,.$$
\label{PB}\end{cor}

By a {\it generalized complex Gaussian measure\/} on $\C^n$, we mean a
generalized Gaussian measure $\ga_\De^c$ on $\C^n\equiv R^{2n}$ with
moments
$$\langle z_j\rangle_\ga =0,\quad \langle z_j z_k\rangle_\ga =0,\quad \langle
z_j\bar z_k\rangle_\ga =\De_{jk},\qquad 1\le j,k\le n,$$ where
$\De=(\De_{jk})$ is an $n\times n$ positive semi-definite hermitian matrix;
i.e. $\ga_\De^c=\ga_{\half \De^c}$,
where $\De^c$ is the $2n\times 2n$ real symmetric matrix of the inner
product on $R^{2n}$ induced by $\De$. As we are interested here in complex
Gaussians, we shall henceforth drop the `$c$' and write
$\ga_\De^c=\ga_\De$.  In particular, if $\De$ is  a strictly positive
hermitian matrix, then
$$\ga_\De= \frac{e^{-\langle\De\inv
z,\bar z\rangle}}{\pi^{n}{\det\De}} d\lcal(z)\,,
$$ where $\lcal$ denotes Lebesgue measure on $\C^n$.

\subsection {Proof of Theorem \ref{usljpd-sphere}}{\label{s-jpd} 
We return
to
our complex Hermitian line bundle $(L,h)$ on a compact almost complex
$2m$-dimensional symplectic manifold
$M$
with symplectic form $\omega=\frac{i}{2}\Theta_L$, where $\Theta_L$ is the
curvature of $L$ with respect to a connection $\nabla$.
We now describe the $n$-point joint distribution arising from our 
probability space $(SH^0_J(M,L^N),\nu_N)$.
Recalling (\ref{inner}), we have the
Hermitian inner product on $H^0_J(M,L^N)$:
\begin{equation*}\label{inner2}\langle s_1, s_2 \rangle = \int_M h^N(s_1,
s_2)\frac{1}{m!}\om^m \quad\quad (s_1, s_2 \in
H^0_J(M,L^N)\,)\;,\end{equation*}
and we write
$\|s\|_2=\langle s,s \rangle^{1/2}$. Recall that $SH^0_J(M,L^N)$ denotes the
unit sphere $\{\|s\|=1\}$ in $H^0_J(M,L^N)$ and $\nu_N$ denotes its Haar
probability measure.

We let $J^1(M,L^N)$ denote the space of
1-jets of sections of
$L^N$. Recall that we have the exact sequence of vector bundles
\begin{equation}\label{jet} 0\to T^*_M\otimes L^N   \stackrel{\iota}{\to}
J^1(M,L^N) \stackrel{\rho}{\to}
L^N\to 0\,.\end{equation} We consider the jet
maps
$$J^1_z:H^0_J(M,L^N)\to J^1(M,V)_z\,,\quad J^1_zs=\ \mbox{the
1-jet
of}\ s\ \mbox{at}\ z\,, \quad \mbox{for}\ z\in M\,.$$ 
The covariant derivative $\nabla:J^1(M,L^N)\to T^*_M\otimes
L^N$  provides a splitting 
of (\ref{jet}) and an  isomorphism
\begin{equation}\label{splitting}(\rho,\nabla):J^1(M,L^N){\buildrel{\approx}
\over
\longrightarrow}  L^N\oplus(T^*_M\otimes L^N)\,.\end{equation}

\begin{defin} 
The {\it $n$-point joint probability distribution\/} at
points $P^1,\dots, P^n$ of $M$ is the probability measure 
\begin{equation}\label{alternately}
\D^N_{(P^1,\dots,P^n)}:=(J^1_{P^1}\oplus\cdots
\oplus J^1_{P^n})_*\nu_N\end{equation} on the space
$J^1(M,L^N)_{P^1}\oplus\cdots \oplus J^1(M,L^N)_{P^n}$.
\end{defin}

Since we are interested in the scaling limit of $\D^N$, we need to describe
this measure more explicitly:  Suppose that
$P^1,\dots,P^n$ lie in a coordinate neighborhood of a point $P_0\in M$ and
choose  preferred coordinates $(z_1,\dots,z_m)$ and a
preferred frame $e_L$ at $P_0$. We let $z^p_1,\dots,z^p_m$ denote the
coordinates of the point $P^p$ ($1\le p\le n$), and we write
$z^p=(z^p_1,\dots,z^p_m)$. (The coordinates of
$P_0$ are $0$.) We consider the
$n(2m+1)$  complex-valued random
variables $x^p,\ \xi^p_q$ ($1\le p\le n,\ 1\le q\le 2m$)
on $S\hcal^2_N(X)\equiv SH^0_J(M,L^N)$ given by 
\begin{equation}\label{dms1def} x^p(s) = s(z^p,0)\,,
\end{equation}  \begin{equation}
\label{dms2def} \xi^p_q(s)=\frac{1}{\sqrtn} \frac{\d^h s}{\d
z_q}(z^p)\,,\quad
\xi^p_{m+q}(s)=\frac{1}{\sqrtn} \frac{\d^h s}{\d\bar z_q}(z^p)\qquad (1\le
q\le m)\,,
\end{equation}
for $s\in SH^0_J(M,L^N)$.

We now write $$x=(x^1,\dots,x^p)\,,\quad
\xi=(\xi^p_q)_{1\le p\le n,1\le q\le 2m}\,,\quad  z=(z^1,\dots,z^n)\,.$$ 
Using (\ref{splitting}) and the variables $x^p,\ \xi^p_q$ to make the
identification
\begin{equation}\label{identification}
J^1(M,L^N)_{P^1}\oplus\cdots \oplus J^1(M,L^N)_{P^n}\equiv
\C^{n(2m+1)}\,,\end{equation} we can write
$$\D^N_z = D^N(x,\xi;z)dxd\xi\,,$$ where
$dxd\xi$ denotes Lebesgue measure on $\C^{n(2m+1)}$.

Before proving Theorem \ref{usljpd-sphere} on the scaling limit of $\D^N_z$,
we state and prove a corresponding result replacing $(SH^0_J(M,L^N),\nu_N)$
with the essentially equivalent Gaussian space $H^0_J(M,L^N)$ with the
{\it normalized standard Gaussian measure}
\begin{equation}\label{gaussian}\mu_N:=\tilde\ga_{2d_N}=
k_{2d_N}e^
{-d_N|c|^2}d\lcal(c)\,,\qquad s=\sum_{j=1}^{d_N}c_jS_j^N\,,\end{equation}
where
$\{S_j^N\}$ is an orthonormal basis for $H^0_J(M,L^N)$.  Recall that this
Gaussian is characterized by the property that the $2d_N$ real variables
$\Re c_j, \Im c_j$ ($j=1,\dots,d_N$) are independent random variables with
mean 0 and variance $1/2d_N$; i.e.,
\begin{equation}\label{normalized}\langle c_j \rangle_{\mu_N}= 0,\quad
\langle c_j c_k\rangle_{\mu_N} = 0,\quad \langle c_j \bar c_k
\rangle_{\mu_N}=
\frac{1}{d_N}\de_{jk}\,.\end{equation} Our normalization guarantees that the
variance of
$\|s\|_2$ is 1: $$\langle
\|s\|^2_2\rangle_{\mu_N}=1\,.$$ We then consider the {\it Gaussian joint
probability distribution\/}
\begin{equation}
\wt \D^N_{(P^1,\dots,P^n)}= \wt D^N(x,\xi;z)dxd\xi =(J^1_{P^1}\oplus\cdots
\oplus J^1_{P^n})_*\mu_N\,.\end{equation}
Since $\mu_N$ is Gaussian and the map
$J^1_{P^1}\oplus\cdots\oplus J^1_{P^n}$
is linear, it follows  that the joint probability
distribution 
is a generalized Gaussian measure of the form
\begin{equation}\label{Dgaussian}D^N(x,\xi;z)dxd\xi=
\ga_{\De^N(z)}\,.\end{equation}
We shall see below that the covariance
matrix $\De^N(z)$ is given in terms of the \szego kernel.  

We have the following alternate form of Theorem \ref{usljpd-sphere}:

\begin{theo}\label{usljpd} Let $L,M,\om$ be as above and let  $\{z_j\}$ be
preferred coordinates centered at a point $P_0\in
M$. Then
$$\wt \D^N_{(z^1/\sqrtn,\dots, z^n/\sqrtn)} \longrightarrow \D^\infty
_{(z^1,\dots,n^n)}$$ where $\D^\infty_{(z^1,\dots,z^n)}$ is the universal
Gaussian measure (supported on the holomorphic 1-jets) of Theorem
\ref{usljpd-sphere}.
\end{theo}

\begin{proof} The covariance matrix 
$\De^N(z)$ in (\ref{Dgaussian}) is a positive semi-definite $n(2m+1)\times
n(2m+1)$ hermitian matrix. If the map $J^1_{z^1}\oplus\cdots\oplus J^1_{z^n}$
is surjective, then $\De^N(z)$ is strictly positive definite and $\wt
D^N(x,\xi;z)$ is a smooth function. On the other hand, if the map is not
surjective, then $\wt D^N(x,\xi;z)$ is a distribution supported on a linear
subspace.  For example, in the integrable holomorphic case, $\wt D^N(x,\xi;z)$
is supported on the holomorphic jets, as follows from the discussion below.

By (\ref{moments}), we have
\begin{eqnarray}&
\Delta^N(z)=\left(
\begin{array}{cc}
A & B \\
B^* & C
\end{array}\right)\,,\nonumber \\
&A=\big( A^{p}_{p'}\big)=
\big \langle x^p \bar x^{p'}\big\rangle_{\mu_N}\,,\quad
B=\big(B^{p}_{p'q'}\big)=
\big\langle x^p \bar \xi^{p'}_{q'}\rangle_{\mu_N}\,,\quad
C=\big(C^{pq}_{p'q'}\big)=
\big\langle  \xi^{p}_{q}\bar \xi^{p'}_{q'}\rangle_{\mu_N}\,,\label{dms3}\\
& p,p'=1,\dots,n, \quad q,q'=1,\dots,2m\,.\nonumber
\end{eqnarray}
(We note that $A,\ B,\ C$ are $n\times n,\ n\times 2mn,\
2mn\times 2mn$ matrices, respectively; $p,q$ index the rows, and
$p',q'$ index the columns.)

We now describe the the entries of the matrix $\De^N$ in terms of the
\szego kernel. 
We have by (\ref{normalized}) and (\ref{dms3}), writing
$s=\sum_{j=1}^{d_N}c_j S^N_j$, 
\begin{equation}\label{dms4}A^{p}_{p'}
=\big \langle x^p \bar x^{p'}\big\rangle_{\mu_N}
= \sum_{j,k=1}^{d_N}\big\langle c_j \bar c_k\big\rangle_{\mu_N}
S_j^N(z^p,0)\overline{S_k^N(z^{p'},0)}=
\frac{1}{d_N}\Pi^{P_0}_N(z^p,0;z^{p'},0)\,.\end{equation}
We need some more notation to describe the matrices $B$ and $C$: Write
$$\nabla_q=\frac{1}{\sqrtn}\frac{\d^h}{\d z_q}\,,\quad
\nabla_{m+q}=\frac{1}{\sqrtn}\frac{\d^h}{\d\bar z_q}\,,\quad 1\le q\le m\,.$$
As in \S \ref{s-notation}, we let $\nabla^1_q$, resp.\ $\nabla^2_q$, denote
the differential operator on $X\times X$ given by applying $\nabla_q$ to the
first, resp.\ second, factor ($1\le q\le 2m$).  By differentiating
(\ref{dms4}), we obtain
\begin{eqnarray}\label{dms5} B^{p}_{p'q'} &=&  \frac{1}{d_N}
\overline{\nabla}^2_{q'}\Pi^{P_0}_N(z^p,0;z^{p'},0)\,,\\
\label{dms6}
C^{pq}_{p'q'} &=& \frac{1}{d_N}
\nabla^1_q\overline{\nabla}^2_{q'}\Pi^{P_0}_N(z^p,0;z^{p'},0)\,.\end{eqnarray}

It follows from (\ref{dms4})--(\ref{dms6})
and
Theorem \ref{neardiag}, recalling (\ref{dhoriz}) and (\ref{RR}), that
\begin{equation}\label{usldelta}\Delta^N(\frac{z}{\sqrtn})\to
\De^\infty(z)= \frac{m!}{c_1(L)^m}\left(
\begin{array}{cc}
A^\infty(z) & B^\infty(z) \\
B^{\infty}(z)^* & C^\infty(z)
\end{array}\right)\end{equation} uniformly, where 
\begin{eqnarray*} A^\infty(z)^p_{p'} &=& \Pi_1^\H(z^p,0;z^{p'},0) \ = \
\frac{1}{\pi^m} e^{\psi_2(z^p,z^{p'})}\,,\\
B^\infty(z)^{p}_{p'q'} &=&\left\{\begin{array}{ll}\frac{1}{\pi^m}
(z_{q'}-w_{q'}) e^{\psi_2(z^p,z^{p'})} \  &\mbox{for}\quad  1\le q\le m\\
0 & \mbox{for}\quad  m+1\le q\le 2m\end{array}\right.\ ,\\
C^\infty(z)^{pq}_{p'q'} &=&\left\{\begin{array}{ll} \frac{1}{\pi^m}
(\bar w_q
-\bar z_q) (z_{q'}-w_{q'}) e^{\psi_2(z^p,z^{p'})} \  &\mbox{for}\quad  1\le
q,q'\le m\\ 0 & \mbox{for}\quad  \max(q,q')\ge m+1\end{array}\right.\ .
\end{eqnarray*}
(Recall that $\psi_2(u,v) =
u \cdot\bar{v} - \half(|u|^2 + |v|^2)$.)
In other words, the coefficients of $\De^\infty(z)$ corresponding to the
anti-holomorphic directions vanish, while the coefficients corresponding to
the holomorphic directions are given by
the \szego kernel $\Pi_1^\H$ for the reduced Heisenberg group and its
covariant derivatives. 

Finally,  we apply Lemma
\ref{continuity} to (\ref{Dgaussian}) and conclude that
$$\wt\D^N_{z/\sqrtn}=\ga_{\De^N(z/\sqrtn)} \to \ga_{\De^\infty(z)}\,.$$
Thus Theorem \ref{usljpd} holds with $\D^\infty_{z}=
\ga_{\De^\infty(z)}$. \end{proof}

\medskip\noindent {\it Proof of Theorem \ref{usljpd-sphere}:\/} The proof is
similar to that of Theorem \ref{usljpd}. This time we define
$$\De^N= \frac{1}{d_N} \jcal_N \jcal_N^*:H^0(M,L^N)\to \C^{n(2m+1)}\,,$$
where $\jcal_N = J^1_{P^1}\oplus\cdots\oplus J^1_{P^n}$ under the
identification (\ref{identification}).  We see immediately that $\De^N$ is
given by (\ref{dms4})--(\ref{dms6}) and the conclusion follows from Lemma
\ref{spherical-vs-gaussian} and (\ref{usldelta}). \qed

\begin{rem} There are other similar ways to define the joint probability
distribution that have the same universal scaling limits.  One of these is to
use the (un-normalized) standard Gaussian measure
$\ga_{2d_N}$ on
$H^0_J(M,L^N)$ in place of the normalized Gaussian $\mu_N$ in Theorem
\ref{usljpd} to obtain joint
densities $D^N_\#(x,\xi;z)=D^N(\frac{x}{N^{m/2}},\frac{\xi}{N^{m/2}};z)$.
Then we would have instead 
$$D^N_\#(N^{m/2}x,N^{m/2}\xi;N^{-1/2}z)dxd\xi\to \ga_{\De^\infty(z)}\,.$$
Another similar result is to let $\la_N$ denote normalized Lebesgue measure on
the unit ball $\{\|s\|\le 1\}$ in $H^0_J(M,L^N)$ and to let
$\wh\D^N_z=\jcal_{N*} \la_N$.  By a similar argument as above, we also have
$\wh\D^N_{z/\sqrtn}\to \ga_{\De^\infty(z)}$.
\end{rem}

\end{document}